\documentclass{gtart}

\input gtoutput
\volumenumber{3}\papernumber{14}\volumeyear{1999}
\pagenumbers{331}{367}\published{14 October 1999}
\proposed{Steve Ferry}
\seconded{Gang Tian, Walter Neumann}
\received{27 March 1999}\revised{30 July 1999}
\accepted{6 October 1999}

\usepackage{amssymb,amscd}
\let\Bbb\mathbb
\let\frak\mathfrak
\renewcommand{\rk}[1]{{\bf #1}\stdspace}

\newtheorem{thm}{Theorem}
\newtheorem{thrm}{Theorem}[section]

\newtheorem{prop}[thrm]{Proposition} 
\newtheorem{cor}[thm]{Corollary} 
\newtheorem{coro}[thrm]{Corollary}

\newcommand{\text}{\mbox}

\begin{document}
        
\title{Examples of Riemannian manifolds with positive\\curvature      
almost everywhere}
\shorttitle{Manifolds with positive curvature almost everywhere}
 
\author{Peter Petersen\\Frederick Wilhelm}

\address{Department of Mathematics, University of California\\
Los Angeles, CA 90095, USA}

\secondaddress{Department of Mathematics,University of California\\
Riverside, CA 92521-0135, USA}

\asciiaddress{Department of Mathematics, University of California\\
Los Angeles, CA 90095, USA\\\\
Department of Mathematics,University of California\\
Riverside, CA 92521-0135, USA}

\email{petersen@math.ucla.edu\\fred@math.ucr.edu}

\begin{abstract} 
We show that the unit tangent bundle of $S^4$ and a real cohomology
$CP^3$ admit Riemannian metrics with positive sectional curvature 
almost everywhere.  These are the only examples so far with positive
curvature almost everywhere that are not also known to admit 
positive curvature.  
\end{abstract}

\asciiabstract{We show that the unit tangent bundle of S^4 and a 
real cohomology CP^3 admit Riemannian metrics with positive
sectional curvature almost everywhere.  These are the only examples so
far with positive curvature almost everywhere that are not also known
to admit positive curvature.}



\maketitlepage

\section{Introduction}

\vglue -0.05in

A manifold is said to have quasi-positive curvature if the curvature is 
nonnegative everywhere and positive at a point.  In analogy with Aubin's    
theorem for manifolds with quasi-positive Ricci curvature one can use the Ricci
flow to show that any manifold with quasi-positive scalar curvature or 
curvature operator can be deformed to have positive curvature.  

By contrast no such result is known for sectional curvature.  In fact,
we do not know whether manifolds with quasi-positive sectional
curvature can be deformed to ones with positive curvature almost
everywhere, nor is it known whether manifolds with positive sectional
curvature almost everywhere can be deformed to have positive
curvature.  What is more, there are only two examples (\cite{GromMey} and
\cite{Esch3}) of manifolds with quasi-positive curvature that are not also
known to admit positive sectional curvature. The example in \cite{Esch3} is
a fake quaternionic flag manifold of dimension $12$. It is not known
if it has positive curvature almost everywhere.  The example in
\cite{GromMey} is on one of Milnor's exotic $7$-spheres.  It was asserted
without proof, in \cite{GromMey}, that the Gromoll--Meyer metric has
positive sectional curvature almost everywhere, but this assertion was
disproven by Mandell, a student of Gromoll (\cite{Man}, cf also \cite{Wil2}).
Thus there is no known example of a manifold with positive sectional
curvature at almost every point that is not also known to admit
positive curvature.  We will rectify this situation here by proving
the following theorem.

\begin{thm}\label{unit tangent bundle}
The unit tangent bundle of $S^4$ admits a  metric with 
positive sectional curvature at almost every point with
the following properties.
\begin{itemize} 
\item[\bf(\romannumeral1)]  
The connected component of the identity of the isometry group is isomorphic to    
$SO(4)$ and contains a free $S^3$--subaction.         
\item[\bf(\romannumeral2)] 
	The set of points where there are $0$ sectional curvatures 
contains totally geodesic flat $2$--tori and is 
the union of two copies of $S^3 \times S^3$ that intersect along a $S^2 \times S^3$.   
\end{itemize}  
\end{thm}

{\bf Remark}\qua
In the course of our proof we will also obtain a precise description
of the set of $0$--curvature planes in the Grassmannian.  
This set is not extremely complicated, but the authors 
have not thought of a description that is succinct enough 
to include in the introduction.  

 \noindent
{\bf Remark}\qua In the sequel to this paper, \cite{Wil2}, the second author shows that 
the metric on the Gromoll--Meyer sphere can be perturbed to one that has 
positive sectional curvature almost everywhere.  In contrast to \cite{Wil2}
the metric we construct here is a perturbation of a metric that has 
zero curvatures at every point.   
\eject

By taking a circle subgroup of the free $S^3$--action in Theorem 
\ref{unit tangent bundle}(\romannumeral1) we get the following.

\begin{cor}\label{6-dim quot}
There is a manifold $M^6$ with the homology of $CP^3$ but not the 
cohomology ring of $CP^3$, that admits a metric with positive sectional
curvature almost everywhere.  
\end{cor}

There are also flat totally geodesic $2$--tori in $M^6$
so as a corollary of Lemma 4.1 in \cite{Stra} 
(cf also Proposition 3 in \cite{Berger3}) we have the following.   

\begin{cor}\label{no 1st order perturbation}  
There are no perturbations of our metrics on the unit 
tangent bundle and $M^6$ whose sectional
curvature is positive to first order. 

	That is, there is no smooth family of metrics $\{ g_t \}_{t\in {\Bbb R}}$
with $g_0$ the metric in Theorem \ref{unit tangent bundle} or  
Corollary \ref{6-dim quot}
so that  
$$
\frac{d}{dt} sec_{g_{t}} (P) |_{t=0} > 0 
$$
for all planes $P$ that satisfy $sec_{g_0} (P) = 0$.   
\end{cor}

\rk{Remark} Although the unit tangent bundle of $S^{4}$ is a
homogeneous space (see below), the metric of Theorem \ref{unit tangent
bundle} is obviously inhomogeneous. What is more if the unit tangent bundle
of $S^{4}$ admits a metric with positive sectional curvature, then it must
be for some inhomogenous metric see \cite{BB}. The space in Corollary 
\ref{6-dim quot} is a biquotient of $Sp(2)$. It follows from 
\cite[Theorem 6]{Onis}  that $M^6$ does not have the homotopy type of 
a homogeneous space.   
\vspace{0.1in}

Before outlining the construction of our metric we  
recall that the $S^3$--bundles over $S^4$
 are classified by ${\Bbb Z} \oplus {\Bbb Z}$ as follows (\cite{Hat}, \cite{Steen}).
The bundle that corresponds to $(n,m) \in {\Bbb Z} \oplus
{\Bbb Z}$ is obtained by gluing two copies of 
${\Bbb R}^4 \times S^3$ together via the diffeomorphism 
$g_{n,m} \co  ( {\Bbb R}^4 \backslash \{ 0\} ) \times S^3
\longrightarrow ( {\Bbb R}^4 \backslash \{ 0\} ) \times S^3$ 
given by 
\addtocounter{thrm}{1}  
\begin{eqnarray}\label{gluing map} 
g_{m,n}(u,v) \longrightarrow 
(\frac{u}{|u|^2}, \frac{u^m v u^n}{|u|^{n+m}}),  
\end{eqnarray}
 where 
we have identified ${\Bbb R}^4$ with ${\Bbb H}$ and
$S^3$ with $\{ v \in {\Bbb H} \ | \ |v| = 1 \}$. We
will call the bundle obtained from $g_{m,n}$ 
``the bundle of type $(m,n)$'', and we will denote it
by $E_{m,n}$.

Translating Theorem 9.5 on page 99 of \cite{Huse} into 
our classification scheme (\ref{gluing map}) shows 
that the unit tangent bundle is of type $(1,1)$. 
We will show it is also the quotient 
of the $S^3$--action on $Sp(2)$ given by 
\begin{eqnarray*}
A_{2,0} ( \  p, \left(\begin{array}{cc}
		a & b \\
		c & d \end{array}\right) \ ) =
 \left(\begin{array}{cc}
		p a     &  p b         \\
		p c     &  p d     \end{array}\right).
\end{eqnarray*}
(It was shown in \cite{Rig} that this quotient is also the total space of the 
bundle of type $(2,0)$, so we will call it $E_{2,0}$.)

The quotient of the biinvariant metric via $A_{2,0}$ is a normal homogeneous space
with nonnegative, but not positive sectional curvature.  To get 
the metric of Theorem \ref{unit tangent bundle} we use 
the method described in \cite{Cheeg} to perturb the biinvariant
of $Sp(2)$ using the commuting $S^3$--actions
\begin{eqnarray*} 
A^u ( \  p_1, \left(\begin{array}{cc}
		a & b \\
		c & d \end{array}\right) \ ) =
 \left(\begin{array}{cc}   
		p_1 a &  p_1 b \\
		c     & d     \end{array}\right)  \\       
A^d ( \  p_2, \left(\begin{array}{cc}  
		a & b \\
		c & d \end{array}\right) \ ) =
 \left(\begin{array}{cc}  
		 a        &      b \\    
		p_2 c     &  p_2 d     \end{array}\right)  \\ 
A^l ( \  q_1, \left(\begin{array}{cc}
		a & b \\
		c & d \end{array}\right) \ ) =
 \left(\begin{array}{cc}
		 a \bar{q}_1   &   b \\
		c   \bar{q}_1  &   d     \end{array}\right)  \\
A^r ( \  p_1, \left(\begin{array}{cc}
		a & b \\
		c & d \end{array}\right) \ ) =
 \left(\begin{array}{cc} 
	            a &      b \bar{q}_2  \\     
		c     &      d  \bar{q}_2    \end{array}\right).
\end{eqnarray*}
We call the new metric on $Sp(2)$, $g_{ \nu_1, \nu_2, l_{1}^u, l_{1}^d }$, and 
will observe in Proposition \ref{S^3 times S^3 times S^3 times S^3}
that  $A_{2,0}$ is by isometries with respect to 
$g_{ \nu_1, \nu_2, l_{1}^u, l_{1}^d }$.
Our metric on the unit tangent bundle is the one induced by the Riemannian
submersion $(Sp(2), g_{ \nu_1, \nu_2, l_{1}^u, l_{1}^d } )
\stackrel{ q_{2,0} }  { \longrightarrow } Sp(2) / A_{2,0} = E_{2,0} $. 
 
In section 1 we review some generalities of Cheeger's method.  
In section 2 we study the symmetries of 
$E_{2,0}$. 
 In section 3 we analyze the infinitesimal 
geometry of the Riemannian submersion $Sp(2) \stackrel{p_{2,1}}{\longrightarrow} S^7$,
given by projection onto the first column.  This will allow us to 
compute the curvature tensor of the metric, $g_{\nu_1, \nu_2}$, obtained 
by perturbing the biinvariant metric via $A^{l}$ and $A^r$.  
In section 4 we compute the $A$--tensor of the Hopf fibration 
$S^7 \longrightarrow S^4$, because it is the key to the geometry of  $g_{\nu_1, \nu_2}$.
In section 5 we specify  the zero curvatures of 
$g_{\nu_1, \nu_2}$, and in section 6 we describe the horizontal space of
$q_{2,0}\co  Sp(2) \longrightarrow E_{2,0}$ with respect to $g_{\nu_1, \nu_2}$
and hence (via results from section 1) with respect to $g_{\nu_1, \nu_2, l_{1}^u, l_{1}^d}$. 
In section 7 we specify the zero curvatures of $E_{2,0}$, first with respect to
$g_{\nu_1, \nu_2}$ and then with respect to $g_{\nu_1, \nu_2, l_{1}^u, l_{1}^d}$,
proving Theorem \ref{unit tangent bundle}.  In section 8 we establish the 
various topological assertions that we made above,
  that $E_{2,0}$ is the total
space of the unit tangent bundle and that 
while the cohomology modules of $E_{2,0}/S^1$ are the 
same as $CP^3$'s the ring structure is different.
Using these computations we will conclude that $E_{2,0}$ and 
$E_{2,0}/S^1$ do not have the homotopy type of any known example of a 
manifold of positive curvature.

We assume that the reader has a working knowledge of O'Neill's 
``fundamental equations of a submersion'' \cite{O'Neill} and the
second author's 
description of the tangent bundle of $Sp(2)$, \cite{Wil1}, we will adopt
results and notation from both papers, in most cases without further 
notice.   It will also be important for the reader to keep the definitions
of the five $S^3$ actions $A^u$, $A^d$, $A^l$, $A^r$ and $A_{2,0}$ straight.  To assist 
we point out that the letters $u$, $d$, $l$ and $r$ stand for ``up'', ``down'', 
``left'' and ``right'' and are meant to indicate the row or column of $Sp(2)$ 
that is acted upon.  

\medskip
\rk{Acknowledgments}
We are grateful to Claude Lebrun, Wolfgang Ziller, and the referee for several
thoughtful and constructive criticisms of the first draft of this paper.  

The first author is supported in part by the NSF.  Support from
National Science Foundation grant DMS-9803258 is gratefully
acknowledged by the second author.

\cl{\it Dedicated to Detlef Gromoll on his sixtieth birthday}

\section{Cheeger's Method}

In \cite{Cheeg} a general method for perturbing the metric on a manifold, $M$, of 
nonnegative sectional curvature is proposed. Various special cases of this
method were first studied in \cite{Berger2} and \cite{BourDesSent}.

If $G$ is a compact group of isometries of $M$, then we let $G$ act on $G
\times M$ by \addtocounter{thrm}{1}  
\begin{eqnarray}  \label{skew action}
g ( p , m) = ( p g^{-1}, g m).
\end{eqnarray}
If we endow $G$ with a biinvariant metric and $G \times M$ with the product
metric, then the quotient of (\ref{skew action}) is a new metric of     
nonnegative sectional curvature on $M$. It was observed in \cite{Cheeg}, that we
may expect the new metric to have fewer $0$ curvatures and symmetries than
the original metric.  

In this section we will describe the effect of certain
Cheeger perturbations on the curvature tensor of $M$.  
Most of the results are special cases of results of \cite{Cheeg}, we have included them 
because they can be described fairly succinctly and are central to all 
of our subsequent computations.
  
The quotient map for the action (\ref{skew
action}) is 
\addtocounter{thrm}{1}   
\begin{eqnarray}  \label{diagnal quotient} 
q_{G \times M} \co  ( p, m) \mapsto pm. 
\end{eqnarray}
The vertical space for $q_{G \times M}$ at $(p, m)$ is 
\begin{eqnarray*}
V_{ q_{G \times M} } = \{ ( - k, k) \ | \ k \in \frak{g} \}
\end{eqnarray*}
where the $-k$ in the first factor stands for the value at $p$ of the
killing field on $G$ given by the circle action \addtocounter{thrm}{1} 
\begin{eqnarray}  \label{right field on G}
(\exp (t k ), p ) \mapsto p \exp( - kt )
\end{eqnarray}
and the $k$ in the second factor is the value of the killing field 
\addtocounter{thrm}{1}   
\begin{eqnarray}  \label{left field on M}
\frac{d}{dt} \exp ( t k) m
\end{eqnarray}
on $M$ at $m$.

Until further notice all Cheeger perturbations under consideration will have the property:
\addtocounter{thrm}{1}
\begin{eqnarray}\label{Simple Cheeger perturbations}      
\mbox{For all $k_1, k_2 \in \frak{g}$ if $\langle (  k_1, 0), (k_2, 0) \rangle = 0$, 
then $\langle (  0, k_1), (0, k_2) \rangle = 0$,} 
\end{eqnarray}
where ${\frak g}$ is the Lie algebra of $G$.  Notice that  
$A^u$, $A^d$, $A^r$ and $A^l$ have this property.  

In this case 
the horizontal space for $q_{G\times M}$ is the direct   
sum 
\addtocounter{thrm}{1}
\begin{eqnarray}\label{Cheeger horizontal}
H_{q_{G\times M}}=\{(\lambda_{2}^{2}k,\lambda_{1}^{2}k)\ |k\in \frak{g}\}
\oplus (\
\{0\}\times H_{O_G}\ ), 
\end{eqnarray}
where $\lambda_{1}=|(-k,0)|$ and $\lambda_{2}=|(0,k)|$ and $H_{O_G}$ is the space that     
is normal to the orbit of $G$. The image of $(\lambda_{2}^{2}k,\lambda_{1}^{2}k)$ under $
dq_{G\times M}$ is 
\addtocounter{thrm}{1}             
\begin{eqnarray}\label{image of vertical under dq}
dq_{G\times M}(\lambda_{2}^{2}k,\lambda_{1}^{2}k)=\frac{d}{dt}p\exp (\lambda_{2}^{2}kt)\exp
(\lambda_{1}^{2}kt)m|_{t=0}=dL_{p,\ast }( \ (\lambda_{1}^{2}+ \lambda_{2}^{2})k\ ). 
\end{eqnarray}
It follows that the effect of Cheeger's perturbation is to keep $H_{O_G}$
perpendicular to the orbits of $G$, to keep the metric restricted to $%
H_{O_G}$ unchanged and to multiply the length of the vector 
$dL_{p,\ast }(\ k \ )$ by the factor 
\addtocounter{thrm}{1} 
\begin{eqnarray}\label{orbit length change}  
\frac{\sqrt{(\lambda_{2}^{2}\lambda_{1})^{2}+(\lambda_{1}^{2}\lambda_{2})^{2}}}{%
(\lambda_{1}^{2}+\lambda_{2}^{2})\lambda_{2}}=\sqrt{\frac{\lambda_{1}^{4}+\lambda_{2}^{2}\lambda_{1}^{2}}{%
( \lambda_{1}^{2} + \lambda_{2}^2 )^2}} =  
\sqrt{ \frac{\lambda_{1}^2 }{ \lambda_{1}^2 + \lambda_{2}^2 } }
\rightarrow 1\text{ as }
\lambda_{1}\rightarrow \infty . 
\end{eqnarray}
If $b$ is a fixed biinvariant metric on $G$ and $l_1$ is a positive real number, then we let $g_{l_1}$ 
denote the metric we obtain on $M$ via the Riemannian submersion 
$q_{G \times M}\co  G \times M \longrightarrow M$ when the metric on the $G$--factor 
in $G \times M$ is $l_{1}^2 b$.  When $G =S^3$, $b$ will always be the unit metric.  
 As pointed out
in (\ref{orbit length change}),  $g_{l_1}$ converges to the original 
metric, $g$, as $l_1 \rightarrow \infty$.  To emphasize this point we will
also denote $g$ by $g_{\infty}$.  
 
Let $TO_G$ denote the tangent distribution in $M$ to the orbits of $G$.
Then any plane $P$ tangent to $M$ can be written as  
\addtocounter{thrm}{1} 
\begin{eqnarray}  \label{tangent to M}
P = span\{ z + k^a, \zeta + k^b \},  
\end{eqnarray}
where $z, \zeta \in H_{O_{G}}$ and $k^a, k^b \in TO_G$.  
We let $\hat{P}$  denote the plane in $T(G\times M)$   
that is horizontal with respect to $q_{G\times M}$ and 
satisfies $dp_2 (\hat{P}) = P$ where $p_2\co  G \times M \longrightarrow M$
denotes the projection onto the second factor.   If     
$\xi \in TM$, then $\hat{\xi} \in T(G\times M)$ has the 
analogous relationship to $\xi$ as $\hat{P}$ has to $P$.  
As pointed out in \cite{Cheeg} we have the following result. 
\begin{prop}\label{Cheeger's curvature condition} 
Let $\lambda_{1}^{a}$ and $\lambda_{2}^{a}$ denote
the lengths of the killing fields on $G$ and $M$ corresponding to $k^{a}$
via the procedures described in (\ref{right field on G}) and (\ref{left
field on M}). 
Let $\lambda_{1}^{b}$ and $\lambda_{2}^{b}$ have the analogous meaning
with respect to $k^{b}$.  
  
\begin{itemize}
\item[\bf(\romannumeral1)]
If the curvature of $P$ is positive with respect to $g_{\infty}$, then 
the curvature of 
\addtocounter{thrm}{1}
\begin{eqnarray}\label{perturbed P}
dq_{G\times M}( \hat{P} )   =   span\{\; z+\frac{ \lambda_{1}^{a  2}  + 
\lambda_{2}^{a 2} }  { \lambda_{1}^{a  2 } }    k^{a},
\;\zeta +   \frac{ \lambda_{1}^{b 2} + \lambda_{2}^{b  2} } { \lambda_{1}^{b  2 } }   k^{b} \; \}
\end{eqnarray}
is positive with respect $g_{l_1}$.  
\item[\bf(\romannumeral2)]
The curvature of $dq_{G\times M}( \hat{P} )$ is positive with respect to $g_{l_1}$ if the 
$A$--tensor,
$A^{q_{G\times M}}$, of $q_{G\times M}$ is nonzero on $\hat{P}$. 
\item[\bf(\romannumeral3)]
If $G = S^3$, then the curvature of $dq_{G\times M}( \hat{P} )$ is positive if 
the projection of $P$ onto $TO_G$ is nondegenerate. 
\item[\bf(\romannumeral4)] 
	If the curvature of $\hat{P}$ is $0$ and $A^{q_{G\times M}}$ vanishes on
$\hat{P}$, then the curvature of $dq_{G\times M}( \hat{P} )$ is $0$.  
\end{itemize}
\end{prop}

\begin{proof}
(\romannumeral2) and (\romannumeral4) are corollaries of O'Neill's horizontal curvature equation.  

To prove (\romannumeral1) notice that
the curvature of $dq_{G\times M}( \hat{P} )$ is positive if the curvature of its
horizontal lift, $\hat{P}$, is positive.  The curvature of $\hat{P}$ is 
positive if its image, $P$, under $dp_2$ has positive curvature, proving (\romannumeral1). 

Let $p_1\co  G \times M \longrightarrow G$ be the projection onto the first factor. 
The curvature of $\hat{P}$ is also 
positive if its image under $dp_1$ is positively curved.  If $G= S^3$, then this 
is the case, provided the image of $\hat{p}$ is nondegenerate, proving (\romannumeral3).\end{proof}

Using (\ref{orbit length change}) we get the following.   
\begin{prop}\label{new horizontal space} 
Let $A_H \co  H \times M \longrightarrow M$ be an action that is  by isometries 
with respect to both $g_{\infty}$ and $g_{l_1}$.  
Let $H_{A_H}$ denote the distribution of vectors that are   
perpendicular to the orbits of $A_H$.
 
$P$ is in $H_{A_H}$ with respect to $g_{\infty}$ if and only if   
$dq_{G\times M}( \hat{P} )$ is in $H_{A_H}$ with respect to 
$g_{l_1}$.
\end{prop}

\begin{proof}
Just combine our description of $g_{l_1}$ with the observation 
that if we square the expression in (\ref{orbit length change})  
we get the reciprocal of the quantity in (\ref{perturbed P}).  
\end{proof}

Ultimately we will be studying Cheeger Perturbations via commuting 
group actions, $A_{G_1}, A_{G_2}$, that individually have property 
(\ref{Simple Cheeger perturbations}).  Generalizing our formulas 
to this situation is a simple matter once we observe the
following result.

\begin{prop}\label{commuting Cheeger Perturbations}
Let $G_1 \times G_2$ act isometrically on $(M, g)$. 
Fix biinvariant metrics $b_1$ and $b_2$, on $G_1$ and $G_2$.  
Let $g_1$ and $g_2$ be the metrics obtained by doing Cheeger 
perturbations of $(M,g)$ with $G_1$ and $G_2$ respectively.
\begin{itemize} 
\item[\bf(\romannumeral1)] 
	$G_2$ acts by isometries on $(M, g_1)$ and $G_1$ acts by 
isometries on $(M, g_2)$.
\item[\bf(\romannumeral2)] 
	Let $g_{1,2}$ denote the metric obtained by doing the Cheeger perturbation
with $G_2$ on $(M, g_1)$ and let $g_{2,1}$ denote the metric obtained by doing the 
Cheeger perturbation with $G_1$ on $(M, g_2)$.  Then 
$$
g_{1,2} = g_{2,1}.
$$     
In fact $g_{1,2}$ coincides with the metric obtained by doing a 
single Cheeger perturbation of $(M,g)$ with $G_1 \times G_2$.
\item[\bf(\romannumeral3)]
 $( k_{u}^1, k_{u}^2, u)$ is horizontal for  
$G_1 \times G_2 \times M \stackrel{ q_{G_1 \times G_2 \times M} }{ \longrightarrow } M$ with respect 
to $b_1 \times b_2 \times g$ if and only if 
$( k_{u}^1, u)$ is horizontal for 
$G_1 \times M \stackrel{ q_{G_1 \times M} }{ \longrightarrow } M$ with respect 
to $b_1 \times g$ and $( k_{u}^2, u)$ is horizontal for 
$G_2 \times M \stackrel{ q_{G_2 \times M} }{ \longrightarrow } M$ with respect 
to $b_2 \times g$.
\end{itemize}
\end{prop}

\begin{proof}
The proof of (\romannumeral1) is a routine exercise in the definitions
which we leave to the reader.

Part (\romannumeral1) gives us a commutative diagram 
\newcommand{\End}{\operatorname{End}}
\begin{eqnarray*}
\begin{CD}
G_1 \times G_2 \times M  @>id \ \times \ q_{G_1 \times M}>>G_2 \times M  \\
@Vid \ \times \ q_{G_2 \times M}VV                 @VVq_{G_2 \times M}V \\
G_1 \times M    @>q_{G_1 \times M}>>  M 
\end{CD}
\end{eqnarray*}
of Riemannian submersions from which (\romannumeral2) readily follows.  

It follows from the diagram that if 
$( \kappa_1, \kappa_2, u)$ is horizontal for $q_{G_1 \times G_2 \times M}$ with respect to   
 $b_1 \times b_2 \times g$, then 
$( \kappa_1, u)$ is horizontal for $q_{G_1  \times M}$ with respect to 
 $b_1  \times g$ and
$(  \kappa_2, u)$ is horizontal for $q_{ G_2 \times M}$ with respect to 
 $ b_2 \times g$.  This proves the ``only if'' part of (\romannumeral3).

The ``if'' part of (\romannumeral3) follows from the ``only if'' 
part and the observation that $H_{  q_{G_1 \times G_2 \times M} } \cap T (G_1 \times G_2)  = 
0$ via a dimension counting argument.
\end{proof}

Rather than changing the metric of 
$E_{2,0}$ directly with a Cheeger perturbation, we will change the 
metric on $Sp(2)$ and then mod out by $A_{2,0}$.     
The constraint to this approach is that the Cheeger perturbations
that we use can not destroy the fact that $A_{2,0}$
 is by isometries.  

Fortunately it was observed in \cite{Cheeg} that if the
metric on the $G$ factor in $G \times M$ is biinvariant, 
then $G$ acts by isometries with respect to $g_{l_1}$.  
Therefore we have the following result.

\begin{prop}\label{S^3 times S^3 times S^3 times S^3}
Let $g_{\nu_1, \nu_2, l_{1}^u, l_{1}^d}$ denote a 
metric obtained from the biinvariant
metric on $Sp(2)$ via Cheeger's method using the
$S^3 \times S^3 \times S^3 \times S^3$--action, $A^{u} \times A^d \times A^l \times A^r$.     

Then $A^{u} \times A^d \times A^l \times A^r$ is by isometries with respect to
 $g_{\nu_1, \nu_2, l_{1}^u, l_{1}^d}$.       
In particular, $A_{2,0}$ is by isometries  with respect to $g_{\nu_1, \nu_2, l_{1}^u, l_{1}^d}$.
\end{prop}
We include a proof
of Proposition \ref{S^3 times S^3 times S^3 times S^3} even though it follows 
from an assertion on page 624 of \cite{Cheeg}.  We do this to establish 
notation that will be used in the sequel, and because the assertion in 
\cite{Cheeg} was not proven.

\begin{proof}
Throughout the paper we will call the tangent spaces to the orbits of 
$A^l$ and $A^r$, $V_1$ and \label{V_1, V_2, H-page}
 $V_2$.  The orthogonal complement of $V_1 \oplus V_2$ 
with respect to the biinvariant metric will be called $H$.  
According to Proposition 2.1 in \cite{Wil1}, $Sp(2)$ is diffeomorphic to 
the pull back of the Hopf fibration $S^7 \stackrel{h}{\longrightarrow} S^4$ 
via $S^7 \stackrel{a \circ h}{\longrightarrow} S^4$, where $a\co  S^4 \longrightarrow S^4$
is the antipodal map and $S^7 \stackrel{h}{\longrightarrow} S^4$ is 
the Hopf fibration that is given by right multiplication by $S^3$.
Moreover, the metric induced on the pull back by the product of two unit $S^7$'s
is biinvariant.  Through out the paper our computations will   
be based on perturbations of the biinvariant metric, $b_{ \frac{1}{ \sqrt{2} } }$,
induced by 
$ S^7(\frac{1}{\sqrt{2}}) \times  S^7(\frac{1}{\sqrt{2}})$, where 
$ S^7(\frac{1}{\sqrt{2}}) $ is the sphere of radius $\frac{1}{\sqrt{2}}$.

Observe that if $k$ is a killing field on $S^3$ whose length is $1$ with respect to
the unit metric, then the corresponding killing field on $Sp(2)$ 
with respect to either $A^l$ or $A^r$ has length $\frac{1}{\sqrt{2}}$ with 
respect to  $b_{ \frac{1}{ \sqrt{2} } }$.
It follows from this that the quantity 
(\ref{orbit length change}) is constant 
when we do a Cheeger perturbation on $b_{ \frac{1}{ \sqrt{2} } }$ 
via either $A^l$ or $A^r$.  Thus the effect of these Cheeger 
perturbations is to scale $V_1$ and $V_2$, and to preserve the splitting
$V_1 \oplus V_2 \oplus H$ and $b_{ \frac{1}{\sqrt{2}} }|_{H}$.  
The amount of the scaling is $< 1$ and converges to $1$ as 
the scale, $l_1$, on the $S^3$--factor in $S^3 \times Sp(2)$ converges 
to $\infty$ and converges to $0$ as $l_1 \rightarrow 0$.  
We will call the resulting scales on $V_1$ and $V_2$, $\nu_1$ and $\nu_2$, 
and call the resulting metric $g_{\nu_1, \nu_2}$.  With this 
convention the biinvariant metric $b_{ \frac{1}{ \sqrt{2} } }$ is 
$g_{ \frac{1}{ \sqrt{2} } , \frac{1}{ \sqrt{2} } }$.

It follows that $g_{\nu_1, \nu_2}$ is the restriction to $Sp(2)$ of the 
product metric $S^{7}_{\nu_1} \times S^{7}_{\nu_2}$ where 
$S^{7}_{\nu}$ denotes the metric obtained from $S^{7} ( \frac{1}{\sqrt{2}} )$
by scaling the fibers of $h$ by $\nu \sqrt{2}$.  Since $A^l$ and $A^r$
are by symmetries of $h$ in each column, they are by isometries on
$S^{7}_{\nu_1} \times S^{7}_{\nu_2}$ and hence also on $(Sp(2), g_{\nu_1, \nu_2})$.

Let $g_{ l_{1}^u, l_{1}^d }$
denote the metric obtained from $b_{ \frac{1}{ \sqrt{2} } }$ via the   
Cheeger perturbation with $A^u\times A^d$ when the metric on the $S^3\times S^3$--factor 
in $S^3 \times S^3 \times Sp(2)$ is $S^3( l_{1}^u) \times S^3(l_{1}^d)$.
An argument similar to the one above, using rows instead of columns, 
shows that $A^u \times A^d$ is by isometries on $(Sp(2), g_{l_{1}^u, l_{1}^d } )$.

Since $A^u \times A^d$ commutes with $A^l \times A^r$, it follows that
$A^u \times A^d$ acts by isometries on $(Sp(2), g_{\nu_1, \nu_2})$.
Doing a Cheeger perturbation with $A^u \times A^d$ on $(Sp(2), g_{\nu_1,\nu_2})$   
produces a metric $g_{\nu_1, \nu_2, l_{1}^u,l_{1}^d}$ which can also be thought of 
as obtained from $Sp(2)$ via a single Cheeger perturbation with 
$A^l \times A^r \times A^u \times A^d$.  Since $A^l \times A^r$
acts by isometries on $(Sp(2), g_{\nu_1, \nu_2})$ and commutes with 
$A^u \times A^d$, $A^l \times A^r$ acts by isometries with
respect to $g_{\nu_1, \nu_2, l_{1}^u,l_{1}^d}$.
       
But $g_{\nu_1, \nu_2, l_{1}^u,l_{1}^d}$ can also be obtained by first perturbing with 
$A^u \times A^d$ and then perturbing with $A^l \times A^r$.      
So repeating the argument of the proceeding paragraph
shows that $A^u \times A^d$ acts by isometries with
respect to $g_{\nu_1, \nu_2, l_{1}^u,l_{1}^d}$.
\end{proof}

It follows that $q_{2,0}\co Sp(2) \longrightarrow E_{2,0} = Sp(2)/A_{2,0}$ is a Riemannian 
submersion with respect to both
$g_{\nu_1, \nu_2}$ and $g_{\nu_1, \nu_2, l_{1}^u,l_{1}^d}$.  We will
abuse notation and call the induced metrics on $E_{2,0}$, 
$g_{\nu_1, \nu_2}$ and $g_{\nu_1, \nu_2, l_{1}^u,l_{1}^d}$.

\section{Symmetries and their effects on $E_{2,0}$}  
Since $A^l \times A^r$ commutes with $A_{2,0}$ and is by isometries    
on $(Sp(2), g_{\nu_1, \nu_2, l_{1}^u, l_{1}^d } )$, it is by isometries 
on $(E_{2,0},  g_{\nu_1, \nu_2, l_{1}^u, l_{1}^d } )$.  However on the 
level of $E_{2,0}$ it has a kernel that at least contains ${\Bbb Z}_2$.  
To see this just observe that the action of $(-1, -1)$ on $Sp(2)$ via   
$A^l \times A^r$ is the same as the action of $-1$ via $A_{2,0}$.  
It turns out that the kernel is exactly ${\Bbb Z}_2$, 
and that $A^l \times A^r$ induces the $SO(4)$--action whose existence was asserted 
in Theorem \ref{unit tangent bundle}(\romannumeral1).

\begin{prop}\label{A^l times A^r on E_2,0}
\begin{itemize}
\item[\bf(\romannumeral1)]
$A^l$ and $A^r$ act freely on $E_{2,0}$.  
\item[\bf(\romannumeral2)]
$E_{2,0}/A^l$ is diffeomorphic 
to $S^4$ and the quotient map 
$p_{2,0}\co  E_{2,0} \longrightarrow E_{2,0}/A^l$ is the bundle of type $(2,0)$.
\item[\bf(\romannumeral3)]
	$A^r$ acts by symmetries of $p_{2,0}$.  The induced map on     
$S^4$ is the join of the standard ${\Bbb Z}_2$--ineffective $S^3$--action
on $S^2$ with the trivial action on $S^1$.  
\item[\bf(\romannumeral4)] The kernel of the action $A^l \times A^r$ on 
$E_{2,0}$ is generated by $(-1,-1)$, so $A^l \times A^r$ induces an 
effective $SO(4)$--action on $E_{2,0}$.  
\end{itemize}
\end{prop}

\begin{proof}
It is easy to see that $A_{2,0} \times A^l$ and $A_{2,0} \times A^r$ 
are free on $Sp(2)$ and hence that  
$A^l$ and $A^r$ are free on $E_{2,0}$, proving (\romannumeral1).   

 It was observed in \cite{GromMey} that
 $Sp(2)/ ( A_{2,0} \times A^l ) = E_{2,0}/A^l$ is diffeomorphic 
to $S^4$ and in \cite{Rig} that $p_{2,0}$ is the bundle of type $(2,0)$. 
Combining these facts proves (\romannumeral2).   

Part (\romannumeral3) is a special case of Proposition 5.5 in \cite{Wil1}. 
  
It follows from (\romannumeral3) that if $(q,p)$ is in the kernel 
of $A^l \times A^r$, then $p = -1$.  On the other hand, $A^l$ 
is the principal $S^3$--action for $E_{2,0} \stackrel{p_{2,0}}{\longrightarrow} S^4$. 
Combining these facts we see that 
the kernel has order $\leq 2$, and we observed above   
that $(-1, -1)$ is in the kernel.
\end{proof}

The $O(2)$ action on $Sp(2)$,   
$A_{O(2)}\co  O(2) \times Sp(2) \longrightarrow S(2)$, 
that is given by
$$
A_{O(2)} \co  ( A, U) \mapsto AU
$$ 
commutes with $A_{2,0}$ and so 
is by symmetries of $q_{2,0}\co  Sp(2) \longrightarrow E_{2,0}$.
Since it also commutes with $A^l \times A^r$ it acts 
by isometries on $(E_{2,0}, g_{\nu_1, \nu_2})$.  It is also 
acts by symmetries of $p_{2,0}$ according to Proposition 5.5 
in \cite{Wil1}.  We shall see, however, that it is not by isometries on
$(Sp(2), g_{\nu_1, \nu_2, l_{1}^u, l_{1}^d })$.  (Note that it 
commutes with neither $A^u$ nor $A^d$.  This is one 
of the central ideas behind the curvature computations in the proof of 
Theorem \ref{unit tangent bundle}.)   

We may nevertheless use $A_{O(2)}$ to find the $0$--curvatures of  
$(E_{2,0}, g_{\nu_1, \nu_2})$.
(We will then use Proposition \ref{Cheeger's curvature condition} to 
see that most of them are not $0$ with respect to 
$g_{\nu_1, \nu_2, l_{1}^u, l_{1}^d }$.) 

A special case of Proposition 5.7 in \cite{Wil1} is the following.
\begin{prop}\label{isometric reduction}
Every point in $E_{2,0}$ has a point in its orbit under 
$A_{SO(2)} \times A^l \times A^r$ that can be represented in $Sp(2)$ 
by a point of the form 
\addtocounter{thrm}{1}
\begin{eqnarray}\label{modulo isometries}  
\left( \; \left(  \begin{array}{cc}  
\cos t   \\
\alpha \sin t 
\end{array} \right) \right.  , \;  \left.  \left( \begin{array}{cc} 
			\alpha \sin t \\
			\cos t   
		\end{array} \right) \; \right), 
\end{eqnarray}
with $t \in [ 0, \frac{\pi}{4}], re(\alpha) = 0$ and $| \alpha | = 1$. 
\end{prop}   
 
We will call points of the form (\ref{modulo isometries}) representative points.
\vspace{.1in}

\noindent
\rk{Notational Convention}
We have seen that $A^r$, $A^l$ and $A_{SO(2)}$ all induce actions on 
$E_{2,0}$, and that $A^r$ and $A_{SO(2)}$ even induce actions on $S^4$.
To simplify the exposition we will make no notational distinction
between these actions and their induced actions.  Thus $A^r$ 
stands for an action on $Sp(2)$, $E_{2,0}$, or $S^4$.  The space that 
is acted on will be clear from the context.

 \section{The curvature tensor of ($Sp(2)$, $g_{\nu_1, \nu_2}$) }

We will study the curvature of  ($Sp(2)$, $g_{\nu_1, \nu_2}$)
by analyzing the geometry of the
riemannian submersion $p_{2,1}\co  Sp(2) \longrightarrow S^7$ given by projection onto
the first column.  

The moral of our story is that the geometry of $p_{2,1}$ resembles the
geometry of $h\co  S^7 \longrightarrow S^4$ very closely. $p_{2,1}$
is a principal bundle with a ``connection metric''. That is the metric is of
the form
$$
\langle X, \ Y\rangle_{t, \omega} =  \langle dp_E (X), \ dp_E (Y)
\rangle_B + t^2 \langle \omega(X),\ \omega(Y)\rangle_G, 
$$
where $G \hookrightarrow E \stackrel{p_E}{\longrightarrow} B$ is a principal
bundle, $\langle \;, \;\rangle_G$ is a biinvariant metric on $G$, $\langle
\;, \;\rangle_B$ is an arbitrary metric on $B$, and $\omega\co  VE
\longrightarrow TE$ is a connection map. In the case of $p_{2,1}$, $G =
S^3 $ and the action is given by $A^r$.  

It will be important for us to understand how the infinitesimal geometry of $%
p_{2,1}$ changes as $\nu _{1}$ and $\nu_{2}$ change. Combining 2.1--2.3 of
\cite{Nash} with a rescaling argument we get the following. 

\begin{prop}
\label{Nash} Let $G\hookrightarrow E\stackrel{p_{E}}{\longrightarrow }B$ be
a principal bundle with a connection metric $\langle \;,\;\rangle
_{t, \omega }$. Let $\nabla ^{t}$, $A^{t}$ and $R^{t}$ denote the
covariant derivative, $A$--tensor and curvature tensor of $\langle
\;,\;\rangle _{t,\omega }$. If the paremeter $t$ is
omitted from the superscript of one of the objects $\nabla ^{t}$, $A^{t}
$, $R^{t}$, or $\langle \;,\;\rangle _{t,\omega }$, then implicity that
parameter has value $1$. 
  
If $e_{1}$, $e_{2}$ and $Z$ are horizontal fields and $\sigma$, $U$ and $%
\tau $ are vertical fields, then:

\begin{itemize}
\item[\bf(\romannumeral1)]   The fibers of $p_{E}$ are totally geodesic.

\item[\bf(\romannumeral2)]  $A_{e_{1}}^{t}e_{2}=A_{e_{1}}e_{2}$, $%
A_{e_{1}}^{t}\sigma =t^{2}A_{e_{1}}\sigma $, $\nabla
_{e_{1}}^{t}e_{2}=\nabla _{e_{1}}e_{2}$, $(\nabla _{e_{1}}^{t}\sigma
)^{v}=(\nabla _{e_{1}}\sigma )^{v}$, $(\nabla _{\sigma }^{t}\tau
)^{v}=(\nabla _{\sigma }\tau )^{v}$,

\item[\bf(\romannumeral3)]  $\langle R^{t}(e_{1},e_{2})e_{2},e_{1}\rangle
_{t}=\langle
R^{B}(dp_{E}(e_{1}),dp_{E}(e_{2}))dp_{E}(e_{2}),dp_{E}(e_{1})\rangle$\nl   
$-\,3t^{2}\Vert A_{e_{1}}e_{2}\Vert ^{2}$,
  
\item[\bf(\romannumeral4)]  $\langle R^{t}(\sigma ,\tau )\tau ,\sigma \rangle
_{t}=t^{2}\langle R(\sigma ,\tau )\tau ,\sigma \rangle $,
 
\item[\bf(\romannumeral5)]  $\langle R^{t}(e_{1},\sigma )\sigma ,e_{1}\rangle
_{t}=t^{4}  \Vert A_{e_{1}}\sigma \Vert ^{2}$,

\item[\bf(\romannumeral6)]  $\langle R^{t}(\sigma ,\tau )U,e_{1}\rangle _{t}=0
$, and

\item[\bf(\romannumeral7)]  $ \langle R^{t}(e_{1},e_{2})Z,\sigma 
\rangle_{t}= t^{2}   \langle R(e_{1},e_{2})Z,\sigma \rangle $.
\end{itemize}
\end{prop}

The next step is to determine the $A$ and $T$ tensors of $p_{2,1}$.

The vertical space of $p_{2,1}$ is 
\begin{eqnarray*}
V_{p_{2,1}} = \{ \ ( 0, \sigma ) \in TSp(2) \ | \ \sigma \in 
V_{h } \} \equiv V_2,
\end{eqnarray*}
where $V_h$ is the vertical space for the Hopf fibration $h\co  S^7 \longrightarrow S^4$ 
that is given by left quaternionic multiplication.   
We will also denote vectors in $V_{p_{2,1}}$ by 
$(0, \tau )$, and we will often abuse notation and write just $\tau$ for $( 0, \tau )$.

\noindent 
\rk{The T--tensor} Given two vector fields, $(0,\sigma _{1})$, $(0,\sigma _{2})$,  
with values in $V_{p_{2,1}}$ we compute 
$$
T_{(0,\sigma _{1})}^{p_{2,1}}(0,\sigma_{2})=
(0,\nabla _{\sigma _{1}}\sigma _{2})^{h}=(0,0) 
$$
since $V_{h}$ is totally geodesic. Therefore \addtocounter{thrm}{1} 
\begin{equation} 
T\equiv 0.  \label{T-tensor of p_1}
\end{equation}

\rk{The A--tensor} 
$H_{p_{2,1}}$ splits as the direct sum 
$H_{p_{2,1}} = V_1 \oplus H$ where $V_1$ and $H$ are as defined on
page \pageref{V_1, V_2, H-page}.    
So any vector
field with values in $H_{p_{2,1}}$ can be written uniquely as $z + w $
where $z$ takes values in $H$ and $w$ takes values in $V_1$. 
Given two such vector fields $z_1 + w_1$ and $z_2 + w_2$ we compute
\addtocounter{thrm}{1} 
\begin{eqnarray}  \label{A-tensor of p_{k,k-1}} 
A^{ p_{2,1} } _{ z_1 + w_1 } z_2 + w_2 = ( \nabla_{ z_1 + w_1 } z_2 + w_2
)^v = ( \nabla_{ z_1 } z_2 )^v =  \nonumber \\
( \ 0, \ ( \ \nabla^{ S^{7}_{ \nu_{2} } } _{ dp_{2}^2(
z_1) } dp_{2}^2 (z_2) \ )^v ) =  
( \  0, \ \frac{1}{2} [ dp_{2}^2(z_1) , dp_{2}^2 (z_2)]^v \ )= \\
( \ 0, \ \nabla^{1}_{ dp_{2}^{2}( z_1) }
dp_{2}^2 (z_2) \ ),   \nonumber 
\end{eqnarray}
where $p_{2}^2 \co  Sp(2) \longrightarrow S^7$ denotes the projection of   
$Sp(2)$ onto its last factor. The second equality is due to the fact that a
vector in $V_2$ is $0$ in the first entry and a vector
in $V_1$ is $0$ in the last entry.

The upshot of (\ref{T-tensor of p_1}) and (\ref{A-tensor of p_{k,k-1}}) is
that the $T$ and $A$ tensors of $p_{2,1}$ are essentially the ``$T$ and $A$
tensors of $h$ in the last factor''. The only difference is that 
$A^{p_{2,1} }$ has a $3$--dimensional kernel, $V_1$, and 
$A^{ h} $ has kernel $= 0$.

This principal also holds for the vertizontal $A$--tensor. If $\sigma =( 0, \sigma )$ is a vector field with values in $V_2$ and 
$z + w$ is a vector field with values in $H_{p_{2,1}} = H \oplus V_1$, then 
\addtocounter{thrm}{1} 
\begin{eqnarray}  \label{vertical A-tensor of p_{k,k-1}}   
A^{p_{2,1}} _{ z + w } \sigma = ( \ 0, \ \nabla^{S^{7}_{ \nu_{ 2} } } _{ d p_{ 2}^2 (z) } \sigma)^h =  
2 \nu_{2 }^2 ( \ 0, \nabla^{ S^{7}( 1 ) }_{ d p_{2}^2 (z) } \sigma \ )^h.
\end{eqnarray}
The only components of the curvature tensor of $Sp(2)$ which are not
mentioned in (\ref{Nash}) are those of the form $\langle R(e_1,\sigma)e_2,
\tau \rangle$ and $\langle R(\sigma,\tau)e_1, e_2 \rangle$. According to
formulas $\{ 2\}$ and $\{ 2^{\prime}\}$ in \cite{O'Neill} these are
\begin{eqnarray*}
\langle R(e_1,\sigma)e_2, \tau \rangle = - \langle (\nabla_\sigma A )_{e_1}
e_2, \tau \rangle - \langle A_{e_1} \sigma , A_{e_2} \tau\rangle
\end{eqnarray*}
and 
\addtocounter{thrm}{1} 
\begin{eqnarray}  \label{mixed curvature}
\langle R(\sigma,\tau)e_1, e_2 \rangle =  \nonumber \\
- \langle (\nabla_\sigma A )_{e_1} e_2, \tau \rangle - \langle A_{e_1}
\sigma , A_{e_2} \tau\rangle + \langle (\nabla_\tau A )_{e_1} e_2, \sigma
\rangle + \langle A_{e_1} \tau , A _{e_2} \sigma\rangle.
\end{eqnarray}
(We use the opposite sign convention for the curvature tensor than O'Neill.)

If we choose $e_1$ and $e_2$ to be basic horizontal fields the first
equation simplifies to 
\addtocounter{thrm}{1} 
\begin{eqnarray}  \label{fundamental equation}
\langle R(e_1,\sigma)e_2, \tau \rangle =  \nonumber \\
- \langle \nabla_\sigma (A _{e_1} e_2),\tau \rangle + \langle
A_{\nabla_\sigma e_1 } e_2, \tau\rangle + \langle A _{e_1} \nabla^{t}_\sigma
e_2, \tau \rangle - \langle A_{e_1} \sigma , A_{e_2} \tau\rangle =  \nonumber
\\
- \langle \nabla_\sigma (A_{e_1} e_2),\tau \rangle - \langle A_{e_2}  
(\nabla_{e_1} \sigma)^h , \tau \rangle + \langle A_{e_1} (\nabla_{e_2}
\sigma)^h, \tau \rangle - \langle A_{e_1} \sigma , A _{e_2} \tau\rangle = \\
- \langle \nabla_\sigma (A_{e_1} e_2),\tau \rangle + \langle A_{e_1} \sigma
, A_{e_2} \tau\rangle - \langle A_{e_2} \sigma, A_{e_1} \tau \rangle -
\langle A_{e_1} \sigma , A_{e_2} \tau\rangle =  \nonumber \\
- \langle \nabla_\sigma (A_{e_1} e_2),\tau \rangle - \langle A_{e_2} \sigma,
A_{e_1} \tau \rangle,  \nonumber
\end{eqnarray}          
where the superscript $^h$ denotes the component in $H$. It appears in
the formula because the orthogonal complement of $H$ in the horizontal space for 
$p_{2,1}$, $V_1$, is
the kernel of $A$.

Similarly (\ref{mixed curvature}) simplifies to \addtocounter{thrm}{1}   
\begin{eqnarray}  \label{fundamental equation 2}
\langle R(\sigma,\tau)e_1, e_2 \rangle =  \nonumber \\    
- \langle \nabla_\sigma (A_{e_1} e_2),\tau \rangle - \langle A_{e_2} \sigma,
A_{e_1} \tau \rangle + \langle \nabla_\tau (A_{e_1} e_2),\sigma \rangle +
\langle A_{e_2} \tau, A_{e_1} \sigma \rangle.  
\end{eqnarray}
According to O'Neill, the curvatures with exactly three horizontal terms are 
\addtocounter{thrm}{1}
\begin{eqnarray}  \label{2 horizontal}
\langle R(e_1, e_2) e_3, \sigma \rangle = - \langle (\nabla_{e_3} A)_{e_1}
e_2, \sigma \rangle =  \nonumber \\ 
- \langle \nabla_{e_3} (A_{e_1} e_2) , \sigma \rangle + \langle A_{ (
\nabla_{e_3} e_1)^h } e_2, \sigma \rangle + \langle A_{ e_1 } (\nabla_{e_3}
e_2)^h, \sigma \rangle.  
\end{eqnarray}

  Through out the rest of this section we will 
let the superscript $^v$ denote the component in $V_1$.  

By (\ref{A-tensor of p_{k,k-1}}) and the fact that the analogous equation
holds for the Hopf fibration we have 
\addtocounter{thrm}{1} 
\begin{eqnarray}  \label{3 double horiz, one vert}
\langle R(e^{h}_1, e^{h}_2) e^{h}_3, \sigma \rangle =  \nonumber \\
- \langle \nabla_{e_{3}^h} (A_{ e_{1}^{h}} e_{2}^h) , \sigma \rangle +   
\langle A_{ ( \nabla_{e^{h}_3} e^{h}_{1})^{ h } } e_{2}^h, \sigma \rangle +
\langle A_{ e_{1}^h } (\nabla_{e^{h}_{3}} e^{h}_2)^{ h }, \sigma \rangle =0.
\end{eqnarray}
Of course 
\addtocounter{thrm}{1} 
\begin{eqnarray} \label{e^v_1, e_2^v, e^3}
\langle R(e^{v}_1, e^{v}_2) e_3, \sigma \rangle = 0  
\end{eqnarray}
since $V_1$ is the kernel of $A$. Using this again we get
\addtocounter{thrm}{1} 
\begin{eqnarray}\label{e^h_1, e_2^v, e^3_v}
\langle R(e^{h}_1, e^{v}_2) e_{3}^v, \sigma \rangle =
 \langle A_{e^{h}_1} ( \ \nabla_{e_{3}^v} e^{v}_2 \ )^h, \sigma \rangle = 0,
\end{eqnarray}
where the last equality is because $( \ \nabla_{e_{3}^v} e^{v}_2 \ )^h = 0$, 
since the orbits of $A^l$ are totally geodesic.

 Also 
\addtocounter{thrm}{1} 
\begin{eqnarray}  \label{bad 3--a prelim}  
\langle R(e^{h}_1, e^{v}_2) e^{h}_3, \sigma \rangle =  \nonumber \\
- \langle \nabla_{e_{3}^h} (A_{ e_{1}^{h}} e_{2}^v) , \sigma \rangle +
\langle A_{ ( \nabla_{e_{3}^h} e^{h}_{1})^{ h } } e_{2}^v, \sigma \rangle +
\langle A_{ e_{1}^h } (\nabla_{e^{h}_{3}} e^{v}_2)^{ h }, \sigma \rangle = \\
\langle A_{ e_{1}^h } (\nabla_{e^{h}_{3}} e^{v}_2)^{ h }, \sigma \rangle = -  
\langle (\nabla_{e^{h}_{3}} e^{v}_2)^h , A_{e_{1}^h} \sigma \rangle. 
\nonumber
\end{eqnarray}  
Keeping in mind that $H$ is the horizontal space for 
$h \circ p_{2,1}$ and that $V_1$ is in the vertical space
for $h \circ p_{2,1}$     
 we
can rewrite the right hand side of (\ref{bad 3--a prelim}) and get  
\addtocounter{thrm}{1} 
\begin{eqnarray}  \label{bad 3--a}    
\langle R(e^{h}_1, e^{v}_2) e^{h}_3, \sigma \rangle = - \langle 
A^{h \circ p_{2,1} }_{e^{h}_{3}} e^{v}_2 , A_{e_{1}^h} \sigma \rangle.
\end{eqnarray}
Using (\ref{bad 3--a}) and the antisymmetry of the curvature tensor we get %
\addtocounter{thrm}{1} 
\begin{eqnarray}  \label{bad 3--b}  
\langle R(e^{v}_1, e^{h}_2) e^{h}_3, \sigma \rangle = \langle 
A^{h\circ p^{2,1} } _{e^{h}_3} e^{v}_{1} , A_{ e_{2}^h} \sigma \rangle.  
\end{eqnarray}
Combining the first Bianchi identity with (\ref{bad 3--a}) 
and (\ref{bad 3--b}) 
we find that         
\addtocounter{thrm}{1}
\begin{eqnarray}  \label{3 horiz, kernal 3--prelim}
\nonumber
\langle R(e_{1}^h, e_{2}^h) e^{v}_3, \sigma \rangle = 
- \langle R( e^{v}_3, e_{1}^h) e_{2}^h, \sigma \rangle - 
\langle R( e^{h}_2, e_{3}^v) e_{1}^h, \sigma \rangle =   \\
- \langle 
A^{h \circ p_{2,1} }_{e^{h}_{2}} e^{v}_3  , A_{e_{1}^h} \sigma \rangle   
+ \langle 
A^{h \circ p_{2,1} }_{e^{h}_{1}} e^{v}_3 , A_{e_{2}^h} \sigma \rangle.
\end{eqnarray}  
It turns out that the sum of the two terms on the right is always $0$ and hence that
\addtocounter{thrm}{1}   
\begin{eqnarray}  \label{3 horiz, kernal 3--prelim-b}
\langle R(e_{1}^h, e_{2}^h) e^{v}_3, \sigma \rangle = 0.  
\end{eqnarray}
This can be seen by direct computation combining Proposition \ref{a-tensor of h} (below)
with section 6 of \cite{Wil1}, but there are many details to check.  A less direct, 
but much quicker proof begins by using the second equation in \ref{Nash}(\romannumeral2) to observe that
the right hand side of (\ref{3 horiz, kernal 3--prelim}) is $4 \nu_{1}^2 \nu_{2}^2$
times the corresponding curvature quantity for the biinvariant metric.  
Therefore it suffices to prove (\ref{3 horiz, kernal 3--prelim-b}) for the biinvariant metric. 
To do this observe that distributions $V_1$ and $V_2$ have the form 
\begin{eqnarray*} 
V_{1,Q} = \left\{ \ L_{Q, *} \left. \left( \begin{array}{cc} 
					\beta &   0   \\
					0    &  0   \end{array}\right) \; \right| \; 
\beta \in {\Bbb H}\; \; \mbox{and} \; \; re(\beta) = 0 \right\}
\end{eqnarray*}  
\begin{eqnarray*} 
V_{2,Q} = \left\{ \ L_{Q, *} \left. \left( \begin{array}{cc} 
					0&   0   \\
					0    &  \beta \end{array}\right) \; \right| \; 
\beta \in {\Bbb H} \; \; \mbox{and} \; \;  re(\beta) = 0   \right\},
\end{eqnarray*}
where $L_{Q, *}$ denotes the differential of left translation by $Q$.

Thus $[ E_{3}^v, \Sigma ] = 0$, where $E_{3}^v$ and $\Sigma$ denote the 
extensions of $ e^{v}_3$ and $\sigma$ to left invariant fields.  
Therefore
\begin{eqnarray*}
\langle R( e_{1}^h, e_{2}^h ) e^{v}_3, \sigma \rangle = - \frac{1}{4} \langle 
[ E_{1}^h, E_{2}^h ],\; [ E_{3}^v, \Sigma ] \rangle = 0, 
\end{eqnarray*}
where $E_{1}^h$ and $E_{2}^h$ denote the extensions of $e_{1}^h$ and $e_{2}^h$
to left invariant fields, proving (\ref{3 horiz, kernal 3--prelim-b}).  

Combining (\ref{e^v_1, e_2^v, e^3}), (\ref{e^h_1, e_2^v, e^3_v}), and
(\ref{3 horiz, kernal 3--prelim-b}) gives us
\addtocounter{thrm}{1}
\begin{eqnarray}  \label{3 horiz, kernal 3}
\langle R(e_{1}, e_{2}) e^{v}_3, \sigma \rangle = 0.
\end{eqnarray}
A review of the results of this section shows that the only difference
between the infinitesimal geometries of $h$ and $p_{2,1}$ is the fact that
for $p_{2,1}$ the terms (\ref{bad 3--a}) and (\ref{bad 3--b}) are not
always $0$ but for $h$ the curvatures terms with exactly three horizontal
vectors are always $0$. This is the root cause of many of the $0$
curvatures in $(Sp(2), g_{\nu_1, \nu_2})$.

\section{The $A$--tensor of $h$}

In section 3 we showed that the $A$--tensor of $p_{2,1}$ is essentially the 
$A$--tensor of the Hopf fibration, $h$, in the last factor. $A^h$ certainly
ought to be very well known, but the authors are not aware of any computation   
of it in the literature. So for completeness we prove the following.

\begin{prop}
\label{a-tensor of h} Let $N$ be a point in $S^{7}$.

\begin{itemize}
\item[\bf(\romannumeral1)]  The vertical space of $h$ at $N$, $V_{h,N}$, is 
$$
\{N\beta \ |\ \beta \in \Bbb{H},re(\beta )=0\}. 
$$

\item[\bf(\romannumeral2)]  If $z$ is in the horizontal space, $H_{h,N}$, for $%
h $ at $N$, then 
\[
A_{z}^{h}N\beta =z\beta . 
\]
\end{itemize}
\end{prop}

\begin{proof}
A proof of this for the case of the complex Hopf fibration can
be found in \cite{O'Neill}. The proof for the quaternionic case is nearly
identical.

(\romannumeral1) is an immediate consequence of the
definition $h$. 

To prove (\romannumeral2) we let a superscript $^h$ denote the horizontal 
component with respect to $h$ and we compute 
\begin{eqnarray*}
( \nabla^{S^7}_{z} N \beta )^h = ( \nabla^{\Bbb{R}^8 }_{z} N \beta )^h = (
(\nabla^{\Bbb{R}^8 }_{z} N) \beta )^h = ( z \beta )^h = z \beta,
\end{eqnarray*}
where for the last equality we have used the fact that $H_{h, N}$ is
invariant under right quaternionic multiplication.
\end{proof}

\section{The Curvature of Generic Planes in Sp(2)}

Because the metric $g_{\nu_1, \nu_2}$ is a simpler metric than 
$g_{\nu_1,\nu_2 l_{1}^u, l_{1}^d}$, we study its $0$ curvatures in this 
section.  Where all statements about curvatures  
of $Sp(2)$ are understood to be with respect to
$g_{\nu_1, \nu_2}$, for some $\nu_1, \nu_2 < \frac{1}{\sqrt{2}}$.

A plane $P$ tangent to $Sp(2)$ has the form 
\addtocounter{thrm}{1}   
\begin{eqnarray}\label{generic form of P}
P = span \{ z + w , \zeta + v \},
\end{eqnarray}
where $z , \zeta \in H$ and $w, v \in V_1 \oplus V_2$. 


\begin{thrm}\label{generic planes} 
\mbox{}\newline
\vspace{-.18in} \mbox{}
\begin{itemize}
\item[\bf(\romannumeral1)]  The curvature of $P$ is positive if the projection 
of $P$ on to one of $H$, $V_{1}$, or $V_{2}$ is two dimensional.

\item[\bf(\romannumeral2)]  If all three projections in (\romannumeral1) are 
degenerate, then the curvature of $P$ is $0$, if its projection onto $H$ is $%
0$.

\item[\bf(\romannumeral3)]  
If all three projections in (\romannumeral1) are
degenerate and the projection of $P$ onto $H$ is one dimensional, then 
we may assume that 
$z\not= 0$ and $\zeta =0$.   The curvature of $P$ is
positive unless   
\addtocounter{thrm}{1} 
\begin{eqnarray}  \label{bad inner product--2}
A^{ h \circ p_{2,1} }_{z} v^1  +
A^{ p_{1,2} }_{ z} v^2  = 0.
\end{eqnarray}

\item[\bf(\romannumeral4)]  If all of the hypotheses of (\romannumeral3) hold  
and (\ref{bad inner product--2}) holds, then the curvature of $P$ is $0$.
\end{itemize}
\end{thrm}   

\noindent
\rk{Remark}
Since $v^2$ is vertical for $h \circ p_{1,2}$ and $A^{ p_{1,2} }_{ z} v^2  = 
A^{ h \circ p_{1,2} }_{ z} v^2$,  
 (\ref{bad inner product--2}) can be thought of 
as $A^{ h \circ p_{1,2} }_{ z} (v^1 + v^2)  = 0$.  Since the fibers of 
$h \circ p_{1,2}$ are totally geodesic this version of 
(\romannumeral3) is a corollary of \ref{Nash}(\romannumeral5).  
We give an alternative proof below.  
\vspace{.1in}

\begin{proof}[Proof of (\romannumeral3) and (\romannumeral4)]
If the projection of $P$ onto $H$, $V_1$ and $V_2$ is degenerate, then
by replacing the second vector with the appropriate linear combination of 
the vectors in (\ref{generic form of P}) we may assume that $P$ 
has the form   
\addtocounter{thrm}{1}
\begin{eqnarray}\label{reduced P}
P = span\{ z + w^1 + w^2, \ v^1 + v^2 \},   
\end{eqnarray}
where $w^1, v^1 \in V_1$ and $w^2, v^2 \in V_2$ are multiples of each other.  

By replacing the first vector in (\ref{reduced P}) by the appropriate 
vector in $P$ we can further assume that $w^2 = 0$, so $P$ has the 
form  
$$
P = span\{ z + w^1   , v^1 + v^2 \},
$$
where $w^1$ and $v^1$ are multiples of each other.  

It follows that 
\addtocounter{thrm}{1}
\begin{eqnarray}\label{vertizontal zeros}
\nonumber
\langle \; R( z + w^1 , \;  v^1 + v^2 ) \;
v^1 + v^2 , \;  z + w^1 \; \rangle = \\
\nonumber
\langle \; R( z , \;  v^1  ) \;  v^1  , \;  z  \; \rangle + 2
\langle \; R( z , \;  v^1  ) \; v^2 , \;  z  \; \rangle + 
\langle \; R( z  , \;  v^2 ) \; v^2 , \;  z  \; \rangle = \\
\| \ A^{h \circ p_{2,1} }_{ z } v^1  \ \|^2 + 
  2 \langle \ A^{ h \circ p_{2,1} }_{z} v^1  , \     
A^{ p_{1,2} }_{ z} v^2 \rangle + 
\| \ A^{ p_{1,2} }_{ z} v^2 \ \|^2 = \\
\nonumber
\| \ A^{ h \circ p_{2,1} }_{z} v^1  +
A^{ p_{1,2} }_{ z} v^2 \ \|^2,
\end{eqnarray}
where we have used the results of section 3 to conclude that 
many terms are $0$.  

(\romannumeral3) and (\romannumeral4) follow from (\ref{vertizontal zeros}).
\end{proof}

\noindent
\rk{Notation}
Since the $A$--tensors in (\ref{vertizontal zeros}) are essentially the 
$A$--tensors of $h$ in the first and second factor we will shorten the notation
and set  
\begin{eqnarray*}
A^1 = A^{ h \circ p_{2,1} } \; \; \mbox{and} \\
A^2 = A^{ p_{1,2} },
\end{eqnarray*}
and we will denote the inner product in (\ref{vertizontal zeros}) by 
\begin{eqnarray*}
\langle \ A^{ h \circ p_{2,1} }_{z} v^1  , \ 
A^{ p_{1,2} }_{ z} v^2\rangle  = 
\langle \ A^{ 1 }_{z} v^1  , \ 
A^{ 2 }_{ z} v^2 \rangle_a.
\end{eqnarray*}
The idea behind this notation is that we are taking a certain type
of inner product between the $A$--tensors of $h$ in the first and second 
factors.  The subscript $_a$ is meant to remind us of the role that the 
antipodal map of $S^4$ plays in determining this inner product.

\rk{Proof of (\romannumeral1)} 
If the projection of $P$ onto either $V_1$ or $V_2$ is nondegenerate,
then according to \ref{Cheeger's curvature condition}(\romannumeral3) 
the curvature of $P$ is positive.   

So we may assume that the projection of $P$ onto $H$ is nondegenerate 
and the projections onto both $V_1$ and $V_2$ are degenerate.  
Since the projection of $P$ onto $V_2$  
is degenerate we may replace the first vector with another vector 
in $P$ that satisfies    
\addtocounter{thrm}{1} 
\begin{eqnarray}  \label{w_2 = 0}
w^2 = 0.
\end{eqnarray}
Since the projection onto $V_1$ is also degenerate 
we may replace our vectors by appropriate linear combinations
to get that either 
\addtocounter{thrm}{1} 
\begin{eqnarray}  \label{v^ = 0 or w^1 = 0}
v^1 = 0 \; \mbox{ or} \; w^1 = 0.
\end{eqnarray}
If $v^1 = 0$, $P$ is spanned by 
$$
\{ z + w^1, \zeta + v^2 \}.  
$$
It follows that 
\begin{eqnarray*}
\langle R  
(z + w^1  , \zeta + v^2  )
\zeta + v^2  , z + w^1     
\rangle_{g_{\nu_1, \nu_2}} \geq \\
 | z|^2 | \zeta^{\perp}|^2 + | A_{\zeta}^{1} w^1 |^2 + 
|A_{z}^2 v^2 |^2 + \\
\langle R(z, \zeta) v^2  , w^1  
 \rangle + \langle R( z , v^2  ) \zeta, 
w^1  
\rangle + \langle R( w^1   , v^2   ) 
\zeta, z \rangle + \langle R( w^1   , \zeta)
v^2  , z \rangle + \\
\langle R( z, v^2 )v^2, w^1 \rangle + \langle R( w^1, v^2 )v^2, z \rangle,
\end{eqnarray*}
where $\zeta^{\perp}$ denotes the component of $\zeta$ that is perpendicular to 
$z$, and we have relied heavily on results of the previous two sections to 
conclude that many of the components of this curvature are $0$.    

Using (\ref{3 horiz, kernal 3}) and the symmetries of the curvature 
tensor
we see that the fourth and the sixth terms on the right are $0$.        
The last two terms are also $0$.  One way to see this is to combine 
(\ref{fundamental equation}) with the fact that $w^1$ is in the 
kernel of $A^{p_{2,1}}$.

Use (\ref{bad 3--a}) to evaluate the other two terms to get 
\addtocounter{thrm}{1}      
\begin{eqnarray}  \label{prelim v^1 =0}
\langle R
(z + w^1  , \zeta + v^2  )       
\zeta + v^2  , z + w^1  
\rangle_{g_{\nu_1, \nu_2}}  \geq  \nonumber
\\
 | z|^2 | \zeta^{\perp}|^2 + | A_{\zeta}^{1} w^1 |^2 + |A_{z}^2 v^2 |^2 - 2 \langle
A_{z}^1 w^1, A_{\zeta} v^2 \rangle_a.
\end{eqnarray}  
For $\xi \in TSp(2)$ let $| \xi |_{1,1} = | dp_{2,1} \xi |_1$ and 
$| \xi |_{1,2} = | dp_{2}^2 \xi|_1$ where $|\cdot|_1$ denotes the unit 
norm on $S^7$.   Combining \ref{Nash}(\romannumeral5), (\ref{vertical A-tensor of p_{k,k-1}}) 
and (\ref{a-tensor of h}) and meticulously chasing through the definitions we get that
\addtocounter{thrm}{1}
\begin{eqnarray}\label{evaluating middle terms} \nonumber
| A_{\zeta}^{1} w^1 |^2 = \nu_{1}^4 \ | \zeta|_{1,1}^{2} \ | w^1 |_{1,1}^2, \\
 |A_{z}^2 v^2 |^2   = \nu_{2}^4 \ | z |_{1,2}^{2} \ | v^2 |_{1,2}^2, \; \; \mbox{and} \hspace*{-.31in} \\
| \ 2 \langle A_{z}^1 w^1, A_{\zeta} v^2 \rangle_a \ | \leq 
 2 \nu_{1}^2 \ \nu_{2}^2 \  | \zeta|_{1,1} \ | w^1 |_{1,1} \ | z |_{1,2} \ | v^2 |_{1,2}. \nonumber
\end{eqnarray}
From (\ref{evaluating middle terms}) we see that sum of the last three terms in
(\ref{prelim v^1 =0}) is always nonnegative so 
\addtocounter{thrm}{1} 
\begin{eqnarray}  \label{final story v^1 =0}
\langle R 
(z + w^1  , \zeta + v^2  )
\zeta + v^2  , z + w^1  
\rangle_{g_{\nu_1, \nu_2}} 
\geq   | z|^2 | \zeta^{\perp} |^2 >0,
\end{eqnarray}    
where the last inequality is due to our hypothesis that the projection of $P$
onto $H$ is nondegenerate.

It remains to consider the case $w^1 = 0$ in (\ref{v^ = 0 or w^1 = 0}). If
so, our plane is spanned by 
$$
\{ z , \zeta + v^1 + v^2 \}.  
$$
Arguing
as in (\ref{prelim v^1 =0}) and (\ref{final story v^1 =0}) we can 
show that 
$$\langle R(z , \zeta + v^1 + v^2   )
\zeta + v^1 + v^2  , z \rangle \geq$$
$$| z|^2 | \zeta^{\perp}|^2 + | A_{z}^{1} v^1 |^2 + |A_{z}^2 v^2 |^2 + 2 \langle
A_{z}^1 v^1, A_{z}^2 v^2 \rangle_a \geq$$
$$| z|^2 | \zeta^{\perp}  |^2 > 0.\eqno{\qed}$$

\begin{proof}[Proof of (\romannumeral2)]  
If $P = span\{ v^1, w^2 \}$ where $v^1 \in V_1$ and $w^2 \in V_2$, then $curv(P)$ $=\,0$ 
because the $T$--tensor of $p_{2,1}$ is $0$ and $V_1$ is the kernel of the 
$A$--tensor of $p_{2,1}$.  
\end{proof}  

\section{Horizontal space of $q_{2,0}$}

From this point forward we will use the notation from section 
6 in \cite{Wil1} for specific tangent vectors to $Sp(2)$.    
Recall that $t \in [0, \frac{\pi}{4}]$ is one of the coordinates 
of our representative points (see Proposition \ref{modulo isometries}).

\begin{prop}
\label{H_p_m,-1} 
For $t \in [0, \frac{\pi}{4})$ the horizontal 
space of $q_{2, 0}$ with respect to $g_{\nu_1, \nu_2}$ 
 is spanned by 
\begin{eqnarray*}
\{ \; 
(x, x), \; \; (- y, y), \; \; ( \eta_{1}, \eta_{1} + \tan(2t) \frac{%
\vartheta_1}{\nu_{2}^2} ) , \; \; (\eta_{2}, \eta_2 + \tan (2t) \frac{
\vartheta_2}{\nu_{2}^2} ) , \\
( - \frac{\frak{v} }{ \nu_{1}^2 }, \; \frac{\frak{v} }{ \nu_{2}^2 }  \; ), \; \; 
( - \frac{\vartheta_{1} }{ \nu_{1}^2 }, \; \frac{\vartheta_{1} }{ \nu_{2}^2 } \; ), \; \; 
( - \frac{\vartheta_{2} } { \nu_{1}^2 }, \; \frac{\vartheta_2}{ \nu_{2}^2 } \; ) \; \}.
\end{eqnarray*}
\end{prop}

\rk{Notation} We will call the seven vectors in (\ref{H_p_m,-1}) $
x^{2, 0}, \; \; y^{2, 0}, \; \; \eta_{1}^{2, 0},$ $\eta_{2}^{2, 0}$, $\frak{v%
}^{2, 0}$, $\vartheta_{1}^{2, 0},$ and $\vartheta_{2}^{2, 0}$ respectively.
We will call the set consisting of the first four $basis(H_{2,0})$, and the
set consisting of the last three $basis(V_{2, 0})$. \vspace{.1in}

\noindent 
\begin{proof}[Sketch of Proof]
The proof is just straight forward computations of inner products. 
Many of the required computations were done in Proposition 6.5 of \cite{Wil1}.
\end{proof} 

A corollary of \ref{generic planes}(\romannumeral1) is the following.
\begin{coro}\label{generic planes in E_2,0}
Let $P$ be a plane in $E_{2,0}$ whose $t$--coordinate is 
in $[0, \frac{\pi}{4})$.  

If the horizontal lift of $P$ with respect to $g_{\nu_1, \nu_2}$  
has a nondegenerate projection onto either 
$H_{2, 0}$ or $V_{2, 0}$, then $P$ is positively curved with 
respect to $g_{\nu_1, \nu_2}$.
\end{coro}

\begin{proof}
We will combine
Theorem \ref{generic planes}(\romannumeral1)
 with the following observations.
\begin{itemize}
\item[\bf(a)]
Any plane in $H_{2,0}$ has a nondegenerate projection 
(with respect to the splitting $H \oplus V_1 \oplus V_2$) 
onto $H$.
\item[\bf(a$'$)]
$V_{2,0} \subset V_1 \oplus  V_2$.
\end{itemize}
It follows from (a) and (a$'$) that if the projection of $P$ onto    
$H_{2,0}$ is nondegenerate, then so is the projection onto   
$H$ (with respect to the splitting $H \oplus V_1 \oplus V_2$). 
\begin{itemize}
\item[\bf(b)] 
Any plane in $V_{2,0}$ has a nondegenerate projection 
(with respect to the splitting $H \oplus V_1 \oplus V_2$)
onto $V_1$.
\item[\bf(b$'$)] 
$ H_{2,0} \subset H \oplus V_2$.   
\end{itemize}   
It follows from (b) and (b$'$) that 
if the projection of $P$ onto $V_{2,0}$ is nondegenerate, then
so is the projection of $P$ onto $V_1$   
(with respect to the splitting $H \oplus V_1 \oplus V_2$). 
\end{proof}
 
\begin{prop} 
\label{H_p_m,-1 , t= pi/4} 
For $t = \frac{\pi}{4}$ the horizontal 
space of $q_{2, 0}$ with respect to $g_{\nu_1, \nu_2}$       
 is spanned by 
\begin{eqnarray*}
(x, x), \; \; (- y, y), \; \; ( \ 0 , \ \frac{\vartheta_1}{\nu_{2}^2} ) , \; \; 
( \ 0, \  \frac{\vartheta_2}{\nu_{2}^2} ) , \\
( - \frac{\frak{v} }{ \nu_{1}^2 }, \; \frac{\frak{v} }{ \nu_{2}^2 }  \; ) , \; \; 
( - \frac{\vartheta_{1} }{ \nu_{1}^2 }, \; \frac{\vartheta_{1} }{ \nu_{2}^2 } \; ) , \; \; 
( - \frac{\vartheta_{2} } { \nu_{1}^2 }, \; \frac{\vartheta_2}{ \nu_{2}^2 } \; ).  
\end{eqnarray*}
\end{prop}    

\begin{proof}[Sketch of proof]
Again the proof is just straight forward computations of inner products.          
In fact the computations that proved (\ref{H_p_m,-1}) will suffice 
for all of the vectors except $( \ 0 , \ \frac{\vartheta_1}{\nu_{2}^2} )$ and  
$( \ 0, \  \frac{\vartheta_2}{\nu_{2}^2} )$.  These are just 
multiples of the limits of $\frac{ \eta^{2,0}_1 } { | \eta^{2,0}_1 | }$  
and $\frac{ \eta^{2,0}_2 } { | \eta^{2,0}_2 | }$ as $t \rightarrow \frac{\pi}{4}$.  
\end{proof} 

\noindent
\rk{Notational Convention}
Let $\pi\co  M \longrightarrow B$ be a riemannian submersion.
To simplify the exposition we will make no notational 
distinction between a vector that is tangent to $B$ and its 
horizontal lift to $M$.  \vspace{.1in}  

\noindent
\rk{Remark}
Since the seven vectors listed in (\ref{H_p_m,-1 , t= pi/4}) are limits 
of (function multiples of) the seven vectors listed in (\ref{H_p_m,-1}) 
we will call them $x^{2, 0}, \; \; y^{2, 0}, \; \; \eta_{1}^{2, 0},$ $\eta_{2}^{2, 0}$, 
$\frak{v}^{2, 0}$, $\vartheta_{1}^{2, 0},$ and $\vartheta_{2}^{2, 0}$.  
Notice that via $q_{2,0} \circ p_{2,0}$, 
$\; x^{2, 0}$ and $y^{2, 0}$ project to the unit normals of 
$S^{2}_{im}$, where $S^{2}_{im}$ is the $2$--sphere in $S^4$ that 
is fixed by the action induced on $S^4$ by $A^r$ 
via  $p_{2,0}\co  E_{2,0} \longrightarrow S^4$ 
(cf the notational remark in section 6 of \cite{Wil1}).  
 Because of this their definitions do not seem to be very canonical 
when $t= \frac{\pi}{4}$.  The approach we have taken is to define them
at representative points with $t= \frac{\pi}{4}$ as extensions of their
definitions at representative points when $t$ is generic.  
We will extend these definitions (as needed) to all points with $t = \frac{\pi}{4}$ 
by letting the isometry group act.   

Also notice that the vectors   
$ (\frac{\vartheta_1}{\nu_{1}^2}, 0 ),(\frac{\vartheta_2}{\nu_{1}^2}, 0 )$ are 
in $H_{q_{2,0}}$ when $t = \frac{\pi}{4}$.  Notice furthermore that at 
$t = \frac{\pi}{4}$
\addtocounter{thrm}{1}      
\begin{eqnarray}\label{extra 0's when t = pi/4}
curv_{\nu_1, \nu_2} ( \; ( \ \frac{\vartheta_i}{\nu_{1}^2}, \  0 \ ),
\eta_{j}^{2,0} \ ) = 0, 
\end{eqnarray}
for any $( \ \frac{\vartheta_i}{\nu_{1}^2}, \  0 \ ) \in span\{
( \ \frac{\vartheta_1}{\nu_{1}^2}, \  0 \ ), ( \ \frac{\vartheta_2}{\nu_{1}^2}, \  0 \ ) \}$
and any $\eta_{j}^{2,0} \in span\{ \eta_{1}^{2,0}, \eta_{2}^{2,0} \}$.    
Here and in the sequel $curv( v,w)$ stands for 
$\langle R(v,w)w,v \rangle$.

\section{Zero Curvatures of $(E_{2,0}, g_{\nu_1, \nu_2})$ and 
$(E_{2,0}, g_{\nu_1, \nu_2, l_{1}^u,l_{1}^d})$ }

It follows from (\ref{generic planes in E_2,0})
 that any $0$--curvature planes for $(E_{2,0}, \ g_{\nu_1, \nu_2})$ have the form  
\addtocounter{thrm}{1}
\begin{eqnarray}\label{difficult planes} 
P = span \{ \zeta, w \},
\end{eqnarray}
where $\zeta \in H_{2,0}$, $w \in V_{2,0}$ and $t \in [0, \frac{\pi}{4})$.         
About such planes we will prove the following.
\begin{prop}\label{0's wrt g_nu_1, nu_2}
\mbox{}\newline
\vspace{-.18in} \mbox{}
\begin{itemize}  
\item[\bf(\romannumeral1)]   
For $t \in [ 0, \frac{\pi}{4})$ the $0$--curvature planes of $(E_{2,0}, g_{\nu_1, \nu_2})$
satisfy (\ref{difficult planes}) and also,  
\addtocounter{thrm}{1}  
\begin{eqnarray}\label{H-part of 0's wrt g_nu_1 , nu_2}    
\zeta \in span\{ x^{2,0}, \eta^{2,0} \}  \nonumber \\ 
w \in span\{ \vartheta^{2,0} \},
\end{eqnarray}
where $\eta^{2,0}$ is a convex combination of 
$\{ \eta_{1}^{2,0}, \eta_{2}^{2,0} \}$ and $\vartheta^{2,0}$
is the same  convex combination of  
$\{ \vartheta_{1}^{2,0}, \vartheta_{2}^{2,0} \}$.  
\item[\bf(\romannumeral2)]
When $t = \frac{\pi}{4}$ the planes of 
the form (\ref{extra 0's when t = pi/4}) have $0$ curvature in $Sp(2)$ with respect to 
$g_{\nu_1, \nu_2}$, and, in addition,
given any convex combination $z$ of $\{ x^{2,0}, y^{2,0} \}$ and any convex combination
$w$ of $\{ \; ( \ - \frac{\vartheta_1}{\nu_{1}^2}, \ 0 \ ), \;  
( \ - \frac{\vartheta_2}{\nu_{1}^2}, \ 0 \ ) \; \}$, 
there is a unique convex combination $v_{w, z}$ of 
$\{ \; ( \ 0 ,\   \frac{\vartheta_1}{\nu_{1}^2} \ ), \; 
( \ 0, \  \frac{\vartheta_2}{\nu_{1}^2} \ ) \; \}$, so that 
\addtocounter{thrm}{1}  
\begin{eqnarray}\label{x-y zeros at pi/4} 
curv_{(Sp(2), g_{\nu_1, \nu_2} )}( z +  \lambda w , w + v_{w,z} ) = 0 
\end{eqnarray}
for all $\lambda \in {\Bbb R}$.  

There are no other zero curvatures when $t =\frac{\pi}{4}$. 
\end{itemize}
\end{prop} 

\begin{proof}
To prove (\romannumeral1) split $\zeta$ into $\zeta = z + \eta$ where $z \in span\{ x^{2,0}, y^{2,0}\}$   
and $\eta \in span\{ \eta_{1}^{2,0}, \eta_{2}^{2,0} \}$.
To be more concrete suppose that $\eta$ is a multiple 
of $\eta_{1}^{2,0}$.  
  
Suppose
$w$ has the form $w = ( - \frac{N_1 \beta}{\nu_{1}^2 }, \frac{N_2 \beta}{\nu_{2}^2}   )$ for some 
purely imaginary unit quaternion $\beta$.  Notice that 
\addtocounter{thrm}{1}
\begin{eqnarray}\label{y inner products}
\langle A_{y^{2,0} }^{1} (- \frac{N_1 \beta}{\nu_{1}^2 }) , \ 
A_{y^{2,0} }^2 \frac{N_2\beta}{\nu_{2}^2}  \rangle =
1, 
\end{eqnarray}
and 
\addtocounter{thrm}{1}
\begin{eqnarray}\label{perp y inner products}
\langle A_{y^{2,0} }^{1} (- \frac{N_1 \beta}{\nu_{1}^2 }) , \ A_{ \xi }^2 \frac{N_2\beta}{\nu_{2}^2}  \rangle =
\langle A_{ \xi  }^{1} (- \frac{N_1 \beta}{\nu_{1}^2 }),  \ 
A_{y^{2,0} }^2 \frac{N_2\beta}{\nu_{2}^2}  \rangle =
0,
\end{eqnarray}
if $\xi \in span\{ x^{2,0}, \eta_{1}^{2,0}, \eta_{2}^{2,0} \}$.   
 Combining \ref{generic planes}(\romannumeral3),(\romannumeral4) with (\ref{y inner products})
and (\ref{perp y inner products}) we get that the curvature of $P$ is positive except possibly if
$\zeta \in span\{ x^{2,0}, \eta \}$.        

On the other hand, the projection of $P$ onto $V_2$ is nondegenerate unless 
$\eta = 0$ or   
\addtocounter{thrm}{1}
\begin{eqnarray}\label{degen proj}  
p_{V_2} (\eta) \; \; \mbox{and $p_{V_2}(w)$ are multiples of each    
other.}
\end{eqnarray}  
So $P$
is positively curved unless $\eta =0$ or $w$ is a multiple of the same  convex combination of 
$\{ \vartheta_{1}^{2,0}, \vartheta_{2}^{2,0} \}$ as $\eta$ is of    
$\{ \eta_{1}^{2,0}, \eta_{2}^{2,0} \}$.  

To complete the proof of (\romannumeral1) it remains to consider the case when $\zeta = x^{2,0}$ 
in which it is easy to see, using 
\ref{generic planes}(\romannumeral3),(\romannumeral4) that $P$ is positively curved 
unless $w \in span \{ \vartheta_{1}^{2,0}, \vartheta_{2}^{2,0} \}$.  

We saw in (\ref{H_p_m,-1 , t= pi/4}) that at points where $t = \frac{\pi}{4}$,
\addtocounter{thrm}{1}
\begin{eqnarray}\label{basis at t= pi/4--b}
\nonumber  
H_{q_{2, 0}} = span \{ 
(x, x), \; \; (- y, y), \; \; ( \ 0 , \ \frac{\vartheta_1}{\nu_{2}^2} ) , \; \; 
( \ 0, \  \frac{\vartheta_2}{\nu_{2}^2} ) , \\
( - \frac{\frak{v} }{ \nu_{1}^2 }, \; \frac{\frak{v} }{ \nu_{2}^2 }  \; ) , \; \;  
( - \frac{\vartheta_{1} }{ \nu_{1}^2 }, 0  \; ) ,  \; \; 
( - \frac{\vartheta_{2} } { \nu_{1}^2 }, \; 0 \; )  \}.
\end{eqnarray} 
Using (\ref{basis at t= pi/4--b}) it is easy to see that the only planes in
$H_{q_{2, 0}} \cap ( V_1 \oplus V_2 )$ whose projections onto both $V_1$ and 
$V_2$ are degenerate are those of the form (\ref{extra 0's when t = pi/4}).

It follows that the only other possible $0$'s are of the type described in 
\ref{generic planes}(\romannumeral3),(\romannumeral4).  Combining this with straightforward computations
of $A$--tensors completes the proof of (\romannumeral2), and even supplies explicit 
solutions for $v_{w,z}$ in terms of $w$ and $z$.  For example when $z = 
x^{2,0}$ 
and $w = ( \ - \frac{\vartheta_1}{\nu_{1}^2}, \ 0 \ )$, $v_{w, z} = 
( \ 0, \ \frac{\vartheta_1}{\nu_{2}^2} \ )$, and when
$z = y^{2,0}$  
and $w = ( \ - \frac{\vartheta_1}{\nu_{1}^2}, \ 0 \ )$, $v_{w, z} = 
( \ 0, \  - \frac{\vartheta_1}{\nu_{2}^2} \ )$.        
\end{proof}

\noindent
\rk{Remark} 
We have shown that the planes in (\ref{0's wrt g_nu_1, nu_2}) have $0$ 
curvature in $(Sp(2),$ $ g_{\nu_1,\nu_2})$.  They also have 
zero curvature in $(E_{2,0}, g_{\nu_1,\nu_2})$.  To verify this     
we have to evaluate the $A$--tensor of $q_{2,0}$ on these planes.  
It turns out that it is $0$.  We will not verify this since 
most of these planes become positively curved with respect to 
$g_{\nu_1,\nu_2, l_{1}^u, l_{1}^d}$.

Now we study the zero curvatures of $(E_{2,0}, g_{\nu_1, \nu_2, l_{1}^u, l_{1}^d})$.       
To do this we adopt the following convention.

\medskip
\noindent
\rk{Notational Convention}
Let $H_{q_{2,0}, \ g_{\nu_1, \nu_2} }$ denote the horizontal space of 
$q_{2,0}$ with respect to 
$g_{\nu_1, \nu_2}$ and $H_{q_{2,0}, \ g_{\nu_1, \nu_2, l_{1}^u, l_{1}^d} }$ the horizontal space of 
$q_{2,0}$ with respect to $g_{\nu_1, \nu_2, l_{1}^u, l_{1}^d}$.  
Let $q_{A^u \times A^d} ^{diag} \co  S^3 \times S^3 \times Sp(2) \longrightarrow Sp(2)$ be the quotient map 
given in (\ref{skew action}) for $A^u \times A^d$.  
According to (\ref{new horizontal space}) $P$ is in $H_{q_{2,0}, \ g_{\nu_1, \nu_2} }$  
if and only if $dq_{A^u \times A^d}^{diag} ( \hat{P})$ is in $H_{q_{2,0} , \ g_{\nu_1, \nu_2, l_{1}^u, l_{1}^d} }$.        
To keep the notation simpler 
we will think of this correspondence as a parameterization of 
$H_{q_{2,0} , \ g_{\nu_1, \nu_2 , l_{1}^u, l_{1}^d } }$
by $H_{q_{2,0}, \ g_{\nu_1, \nu_2} }$ and we will denote vectors and planes in
 $H_{q_{2,0} , \ g_{\nu_1, \nu_2 l_{1}^u, l_{1}^d } }$
by the corresponding vectors and planes in $H_{q_{2,0},\ g_{\nu_1, \nu_2} }$. We 
will do this without any further mention or change in notation.  
For example $x^{2,0}$ now stands for $dq_{A^u \times A^d}^{diag} ( \hat{x}^{2,0})$.  

With this convention it follows from 
\ref{Cheeger's curvature condition}(\romannumeral1) 
that the zero curvatures of $g_{\nu_1, \nu_2, l_{1}^u, l_{1}^d}$
are a subset of the zero curvatures of $g_{\nu_1, \nu_2}$.
We will show that most of the planes in (\ref{0's wrt g_nu_1, nu_2}) have nondegenerate 
projections onto the orbits of $A^u$ and hence are positively curved with 
respect to $g_{\nu_1, \nu_2, l_{1}^u, l_{1}^d }$ according to 
\ref{Cheeger's curvature condition}(\romannumeral3). 
The root cause of this is that $A_{SO(2)}$ is not by isometries with respect to 
$g_{\nu_1, \nu_2, l_{1}^u, l_{1}^d}$.  This seems extremely plausible 
because $A_{SO(2)}$ commutes with neither $A^u$ nor $A^d$.   
A rigorous proof that it does not act by isometries follows from our curvature 
computations.

In any case we no longer know that $A_{SO(2)}$ is by isometries, so (\ref{modulo isometries})
are no longer our representative points.

The representative points are now 
\addtocounter{thrm}{1}\small
\begin{eqnarray}\label{theta,t representatives}
\nonumber
\left(\begin{array}{cc} 
\cos \theta    &    \sin \theta     \\ 
- \sin \theta  &   \cos \theta    \end{array}\right) 
\left( \; \left(\begin{array}{c}
\cos t   \\
\alpha \sin t \end{array} \right) \right.   , 
\left. \; \left(\begin{array}{c}
\alpha \sin t   \\
\cos t \end{array} \right) \; \right) =  \\
\left( \; \left(\begin{array}{c}
\cos \theta \cos t + \alpha \sin \theta \sin t \\
- \sin \theta \cos t + \alpha \cos \theta \sin t  
\end{array}\right) \right.  , \;  
\left.  \left(\begin{array}{c}
\alpha \cos \theta \sin t + \sin \theta \cos t \\
- \alpha \sin \theta \sin t +  \cos \theta \cos t 
\end{array}\right) \right)  \equiv \nonumber
 \\
( \ N_{1}( \theta, t), \ N_{2}(\theta, t) \ ),
\end{eqnarray}\normalsize 
for all $t \in [0, \frac{\pi}{4}]$ and for all $\theta \in [0, \pi)$.  
The interval for $\theta$ is $[0, \pi)$ because $A_{SO(2)}$ is 
${\Bbb Z}_2$--ineffective.     

 Using this we will show the following.

\begin{prop}\label{proj onto  eta}
For any $t \in [0, \frac{\pi}{4})$
the planes of the form 
\addtocounter{thrm}{1}
\begin{eqnarray}\label{generic eta-vartheta 0's}
P = span\{ \eta^{2, 0}, \vartheta^{2,0} \}
\end{eqnarray}
have positive curvature with respect to $g_{\nu_1, \nu_2, l_{1}^u, l_{1}^d}$ 
 unless $\theta = 0, \frac{\pi}{4}, \frac{3\pi}{4},$ or $\frac{\pi}{2}$. 
 
Here $\eta^{2, 0}$ and $\vartheta^{2,0}$ stand for (the same)
convex combinations of $\{ \eta_{1}^{2,0}, \eta_{2}^{2,0}\}$  
and $\{ \vartheta_{1}^{2,0}, \vartheta_{2}^{2,0}\}$
respectively.
\end{prop}

\begin{proof}
The main idea is that most planes are eliminated because their projection 
onto the orbits of $A^u$ are nondegenerate.    
Combining \ref{Cheeger's curvature condition}(\romannumeral3) with 
(\ref{commuting Cheeger Perturbations}) we see that 
it is actually enough to check that a projection is degenerate with 
respect to $b_{ \frac{1}{ \sqrt{2} } }$.

The tangent space to the orbit of $A^u$ is spanned by the 
$b_{ \frac{1}{ \sqrt{2} } }$--orthogonal basis 
\addtocounter{thrm}{1}\small
\begin{eqnarray}\label{basis for TO_A^u}
\nonumber
\{ \; U_{\alpha} , \; U_{\gamma_1}, \; U_{\gamma_2} \; \} \equiv \\
\nonumber
\left\{ \; \; 
\left( \; \left(\begin{array}{c}  
\alpha ( \ \cos \theta \cos t + \alpha \sin \theta \sin t  \ )  \\ 
 0 
\end{array}\right) \right. \right.  , \;  
\left.  \left(\begin{array}{c}  
\alpha ( \alpha \cos \theta \sin t + \sin \theta \cos t \ ) \\
 0 
\end{array}\right) \right),    \\
\nonumber
\left( \; \left(\begin{array}{c}
 \gamma_1 ( \  
\cos \theta \cos t + \alpha \sin \theta \sin t \ ) \\
 0 
\end{array}\right) \right.  , \;  
\left.  \left(\begin{array}{c}
 \gamma_1 ( \ 
\alpha \cos \theta \sin t + \sin \theta \cos t \ )   \\
 0
\end{array}\right) \right),    \\
\left( \; \left(\begin{array}{c}
 \gamma_2 ( \ 
\cos \theta \cos t + \alpha \sin \theta \sin t \ ) \\
 0 
\end{array}\right) \right.  , \;  
\left. \left.  \left(\begin{array}{c}
 \gamma_2
( \ \alpha \cos \theta \sin t + \sin \theta \cos t  \ ) \\
 0
\end{array}\right) \right) \; \; \right\} , \quad\mbox{}
\end{eqnarray}\normalsize
where $\gamma_1$ and $\gamma_2$ are purely imaginary, unit quaternions that satisfy 
$\gamma_1 \gamma_2 = \alpha$.

Before we can find the projections we also need a formula for 
$( \eta_{1}, \eta_1 )$ along the orbits of 
$A_{SO(2)}$.  
At points of the form (\ref{theta,t representatives}) $(\eta_1, \eta_1)$ is 
\small\begin{eqnarray*}
\left(\begin{array}{cc}
\cos \theta    &    \sin \theta     \\
- \sin \theta  &   \cos \theta    \end{array}\right) 
\left( \; \left(\begin{array}{c}
- \gamma_1 \sin t   \\
\gamma_2 \cos t \end{array} \right) \right.  , 
\left. \; \left(\begin{array}{cc}
\gamma_2 \cos t   \\
- \gamma_1 \sin t \end{array} \right) \; \right) =  \\
\left( \; \left(\begin{array}{c}    
-\gamma_1 \cos \theta \sin t + \gamma_2 \sin \theta \cos t \\
\gamma_1 \sin \theta \sin t + \gamma_2 \cos \theta \cos t  
\end{array}\right) \right.  , \; 
\left.  \left(\begin{array}{c}
\gamma_2 \cos \theta \cos t - \gamma_1 \sin \theta \sin t \\
- \gamma_2 \sin \theta \cos t -  \gamma_1 \cos \theta \sin t 
\end{array}\right) \; \right).
\end{eqnarray*}\normalsize
To find the projections we compute 
\addtocounter{thrm}{1}
\small\begin{eqnarray}\label{<(eta, eta), U_gamma_1>}\nonumber
\langle ( \eta_{1}, \eta_{1}) , U_{\gamma_1} \rangle = \\ \nonumber
\left\langle
\left( \; \left(\begin{array}{c}
-\gamma_1 \cos \theta \sin t + \gamma_2 \sin \theta \cos t \\
\gamma_1 \sin \theta \sin t + \gamma_2 \cos \theta \cos t 
\end{array}\right) \right. \right.  , \; 
\left.  \left(\begin{array}{c}
\gamma_2 \cos \theta \cos t - \gamma_1 \sin \theta \sin t \\ 
- \gamma_2 \sin \theta \cos t -  \gamma_1 \cos \theta \sin t  
\end{array}\right) \; \right)  , \\ \nonumber 
\left.
\left( \; \left(\begin{array}{c}
\gamma_1 \cos \theta \cos t - \gamma_2 \sin \theta \sin t \\
0
\end{array}\right) \right. \right.  , \; 
\left.  \left.  \left(\begin{array}{c}
- \gamma_2 \cos \theta \sin t + \gamma_1 \sin \theta \cos t \\   
0 
\end{array}\right) \right) \right\rangle   =  \\
\nonumber \frac{1}{2} ( \ 
- \cos^2 \theta \sin t \cos t - \sin^2 \theta \cos t \sin t 
- \cos^2 \theta \sin t \cos t - \sin^2 \theta \cos t \sin t \ ) =\\
 -  \sin t \cos t = - \frac{1}{2} \sin 2t.\qquad\mbox{}  
\end{eqnarray}\normalsize
According to (\ref{H_p_m,-1}) we should compute the inner  
products of 
$ ( \ - \frac{\vartheta_1}{.5}, \  \frac{\vartheta_1}{.5} \ )$ and 
$( \ - \frac{\vartheta_2}{.5}, \  \frac{\vartheta_2}{.5} \ )$
with the vectors in (\ref{basis for TO_A^u}) with respect to 
$b_{\frac{1}{\sqrt{2}}}$.  These are the same as the inner 
products of  
$( \ - \vartheta_1, \  \vartheta_1 \ )$ and 
$( \ - \vartheta_2, \  \vartheta_2 \ )$ with the vectors in (\ref{basis for TO_A^u}) with respect to 
$b_{1}$, where $b_1 = 2 b_{\frac{1}{\sqrt{2}}}$.  We will compute the 
latter inner products since the notation is simpler.  We will only do
this explicitly for $( \ - \vartheta_1, \  \vartheta_1 \ )$,
since the computations for $( \ - \vartheta_2, \  \vartheta_2 \ )$
are the same modulo obvious changes in notation.  
\small
\addtocounter{thrm}{1}
 \begin{eqnarray}\label{<(0, vartheta), U_gamma_1>} 
\nonumber
\langle ( 0, \vartheta_{1} ) , U_{\gamma_1} \rangle = \\
  \left\langle 
\left( \; \left(\begin{array}{c}
0 \\
0
\end{array}\right) \right. \right.  , \; 
\left.  \left(\begin{array}{c}
\gamma_2 \cos \theta \sin t + \gamma_1 \sin \theta \cos t \\
Blah 
\end{array}\right) \; \right)  , \nonumber  \\  
\nonumber
\left.
\left( \; \left(\begin{array}{c}
\gamma_1 \cos \theta \cos t - \gamma_2 \sin \theta \sin t \\
0  
\end{array}\right) \right. \right.  , \; 
\left.  \left.  \left(\begin{array}{c}
- \gamma_2 \cos \theta \sin t + \gamma_1 \sin \theta \cos t \\
0
\end{array}\right) \right) \right\rangle   =  
\nonumber  \\
- \cos^2 \theta \sin^2 t + \sin^2 \theta \cos^2 t,\qquad\mbox{} 
\end{eqnarray}\normalsize
where ``$Blah$'' stands for a nonzero, but irrelevant term.  

Combining the previous 2 equations we get\small 
\addtocounter{thrm}{1}
\begin{eqnarray}\label{proj eta_1^2-1 onto U_gamma_1}
\langle \eta_{1}^{2,0}, U_{\gamma_1} \rangle = 
- \frac{1}{2} \sin 2t + \tan(2t) ( \ - \cos^2 \theta \sin^2 t + \sin^2 \theta \cos^2 t \ )  
\end{eqnarray}\vspace{-0.1in}
\addtocounter{thrm}{1}
\begin{eqnarray}\label{ < U_gamma_1, (vartheta_1, 0) >} \nonumber
\langle U_{\gamma_1} , ( - \vartheta_1 , 0 ) \rangle = \\  \nonumber
\left\langle \; \;  
\left( \; \left(\begin{array}{c}
\gamma_1 \cos \theta \cos t - \gamma_2 \sin \theta \sin t \\
0
\end{array}\right) \right. \right. , \; 
\left.  \left.  \left(\begin{array}{c}
- \gamma_2 \cos \theta \sin t + \gamma_1 \sin \theta \cos t \\   
0
\end{array}\right) \right) \right. ,  \\
\left. 
\left( \; \left(\begin{array}{c}
-  \gamma_1 \cos \theta \cos t - \gamma_2 \sin \theta \sin t \\
Blah 
\end{array}\right) \right. \right. , \; 
\left.  \left.  \left(\begin{array}{c}
0   \\   
0
\end{array}\right) \right) \right\rangle  =\nonumber  \\  
 - \cos^2 \theta \cos^2 t + \sin^2 \theta \sin^2 t.\qquad\mbox{}  
\end{eqnarray}\normalsize  
Combining (\ref{<(0, vartheta), U_gamma_1>}) and (\ref{ < U_gamma_1, (vartheta_1, 0) >})
we get\small 
\addtocounter{thrm}{1}
\begin{eqnarray}\label{< (- vartheta, vartheta), U_gamma_1 >} 
\langle (- \vartheta_1, \ \vartheta_1), \; U_{\gamma_1} \rangle = \nonumber \\ 
- \cos^2 \theta \sin^2 t + \sin^2 \theta \cos^2 t - \cos^2 \theta \cos^2 t
+ \sin^2 \theta \sin^2 t =  \nonumber  \\
- \cos^2 \theta ( \ \sin^2 t + \cos^2 t \ ) + \sin^2 \theta ( \ \cos^2 t + \sin^2 t \ ) =\nonumber\\
- \cos 2\theta.
\end{eqnarray}
\addtocounter{thrm}{1} \vspace{-0.1in}
\begin{eqnarray}\label{ < U_gamma_2, eta_1 >}
\nonumber
\langle U_{\gamma_2}, \; ( \ \eta_1 , \eta_1 \ ) \rangle = \\ \nonumber
\left\langle \; \;  \left( \; \left( \begin{array}{c}  
\gamma_2 \cos \theta \cos t + \gamma_1 \sin \theta \sin t \\
					0   \end{array}\right) \right. \right. \ , \; 
\left( \left.  \begin{array}{c}
\gamma_1 \cos \theta \sin t + \gamma_2 \sin \theta \cos t \\
0 
\end{array}\right) \; \right),  \\ \nonumber
\left( \; \left(\begin{array}{c}
-\gamma_1 \cos \theta \sin t + \gamma_2 \sin \theta \cos t \\
\gamma_1 \sin \theta \sin t + \gamma_2 \cos \theta \cos t 
\end{array}\right) \right.  , \; 
\left.  \left. \left(\begin{array}{c}
\gamma_2 \cos \theta \cos t - \gamma_1 \sin \theta \sin t \\
- \gamma_2 \sin \theta \cos t -  \gamma_1 \cos \theta \sin t 
\end{array}\right) \; \right)  \; \; \right\rangle = \\
\nonumber\frac{1}{2} ( \  
\cos^2 t \cos \theta \sin \theta - \sin^2 t \cos \theta \sin \theta - 
\sin^2 t \cos \theta \sin \theta + \cos^2 t \cos \theta \sin \theta \ ) 
 =  \\
\frac{1}{4} \sin 2 \theta ( \ 2 \cos^2 t - 2 \sin^2 t \ ) = \nonumber \\
\frac{1}{2} \sin 2 \theta \cos 2 t.\qquad\mbox{}
\end{eqnarray} 
\addtocounter{thrm}{1}\vspace{-0.1in}  
\begin{eqnarray}\label{ < U_gamma_2, vartheeta_1 >}
\nonumber     
\langle U_{\gamma_2}, ( \ 0, \ \vartheta_1 \ ) \rangle = \\ \nonumber
\left\langle \; \;  \left( \; \left( \begin{array}{c} 
\gamma_2 \cos \theta \cos t + \gamma_1 \sin \theta \sin t \\
					0   \end{array}\right) \right. \right. \ , \; 
\left( \left.  \begin{array}{c}
\gamma_1 \cos \theta \sin t + \gamma_2 \sin \theta \cos t \\
0 
\end{array}\right) \; \right),  \\ 
\left. \; \;  \left( \; \left( \begin{array}{c}  
0 \\ 
0   \end{array}\right) \right. \right. \ , \; 
\left( \left. \left.  \begin{array}{c}
\gamma_2 \cos \theta \sin t + \gamma_1 \sin \theta \cos t \\  
blah
\end{array}\right) \; \right)  \; \; \right\rangle  = \nonumber  \\
2 \sin \theta \cos \theta \sin t \cos t  = \frac{1}{2} \sin 2 \theta \sin 2t.
\qquad\mbox{}
\end{eqnarray}\normalsize
Combining (\ref{ < U_gamma_2, eta_1 >}) and (\ref{ < U_gamma_2, vartheeta_1 >}) 
we get\small 
\addtocounter{thrm}{1}
\begin{eqnarray}\label{< U_gamma_2, eta_1^2,0>}
\nonumber
\langle U_{\gamma_2}, \  \eta_{1}^{2,0} \rangle = \\
\frac{1}{2} \sin 2\theta \cos 2t + \frac{1}{2} \tan 2 t \sin 2\theta \sin 2t = \nonumber \\ 
\frac{\sin 2 \theta}{2} ( \  \cos 2t +  \tan 2 t \sin 2t \ ).  
\end{eqnarray}\vspace{-0.1in}
\addtocounter{thrm}{1}  
\begin{eqnarray}\label{ < U_gamma_2, (vartheta_1, 0) >} \nonumber 
\langle U_{\gamma_2} , ( - \vartheta_1 , 0 ) \rangle = \\  \nonumber
\left\langle \; \;  
\left( \; \left(\begin{array}{c}   
\gamma_2 \cos \theta \cos t + \gamma_1 \sin \theta \sin t \\
0 
\end{array}\right) \right. \right. , \; 
\left.  \left.  \left(\begin{array}{c}
\gamma_1 \cos \theta \sin t + \gamma_2 \sin \theta \cos t \\   
0
\end{array}\right) \right) \right. , \\
\left. 
\left( \; \left(\begin{array}{c} 
-  \gamma_1 \cos \theta \cos t - \gamma_2 \sin \theta \sin t \\ 
Blah  
\end{array}\right) \right. \right. , \; 
\left.  \left.  \left(\begin{array}{c}
0   \\   
0 
\end{array}\right) \right) \right\rangle  =   \nonumber  \\ 
-  2 \cos \theta  \sin \theta \cos t \sin t =  - \frac{1}{2} \sin 2 \theta \sin 2 t.\qquad\mbox{} 
\end{eqnarray}\normalsize 
 Combining (\ref{ < U_gamma_2, vartheeta_1 >}) and (\ref{ < U_gamma_2, (vartheta_1, 0) >})
we see that 
\addtocounter{thrm}{1}
\begin{eqnarray}\label{< ( -vartheta_1, vartheta_1), U_gamma_2}
\langle ( -\vartheta_1, \ \vartheta_1 \ ), \ U_{\gamma_2} \rangle = 0.
\end{eqnarray}  
Combining (\ref{< (- vartheta, vartheta), U_gamma_1 >}), (\ref{< U_gamma_2, eta_1^2,0>}), and
(\ref{< ( -vartheta_1, vartheta_1), U_gamma_2}) we see that the projection of 
$P$ onto $TO_{A^u}$ is nondegenerate unless $\cos 2 \theta = 0$ or 
$\sin 2 \theta = 0$, proving the proposition. 
\end{proof}

\noindent
\rk{Remark}
One might hope to show that the planes from (\ref{proj onto  eta}) with $\theta = 0, \frac{\pi}{4},
 \frac{\pi}{2}$, and $\frac{3\pi}{4}$ are positively curved by studying the projections
onto $TO_{A^d}$.  These are also nondegenerate for most values of $\theta$.  Unfortunately
they are degenerate precisely when $\theta = 0, \frac{\pi}{4},
 \frac{\pi}{2}$, and $\frac{3\pi}{4}$.

\begin{prop}\label{x-y, planes} 
Let $\eta^{2,0}$ and $\vartheta^{2,0}$ be as in 
(\ref{proj onto  eta}).  
For $t \in [0, \frac{\pi}{4})$ if
\addtocounter{thrm}{1}
\begin{eqnarray}\label{generic x, vartheta 0's}
P = span\{ \zeta, \vartheta^{2,0} \}
\end{eqnarray}
where $\zeta \in span\{ x^{2,0}, \eta^{2,0} \}$, then the curvature of $P$ is positive
with respect to $g_{\nu_1, \nu_2, l_{1}^u, l_{1}^d}$   
unless $\theta = 0, \frac{\pi}{4}, \frac{\pi}{2},$ or $\frac{3 \pi}{4}$. 
\end{prop}

\begin{proof}
Along the orbit of $SO(2)$ the formula for $x^{2,0}$ is\small 
\begin{eqnarray*}
\left(\begin{array}{cc}
\cos \theta    &    \sin \theta   \\ 
- \sin \theta  &    \cos \theta   
\end{array}\right) 
\left( \; \left(\begin{array}{c} 
	- \sin t    \\ 
\alpha \cos t    \end{array} \right)  \right. , \;  
\left.   \left( \begin{array}{c}   
	\alpha \cos t   \\
	- \sin t   \end{array} \right) \; \right) =  \\    
\left( \; \left(\begin{array}{c} 
- \cos \theta \sin t  + \alpha \sin \theta \cos t \\      
\sin \theta \sin t + \alpha \cos \theta \cos t
\end{array} \right) \right .   , \; 
\left. \left( \begin{array}{c} 
\alpha \cos \theta \cos t - \sin \theta \sin t    \\
- \alpha \sin \theta \cos t - \cos \theta \sin t 
\end{array}\right) \; \right).        
\end{eqnarray*}\normalsize  
So 
\addtocounter{thrm}{1}\small
\begin{eqnarray}\label{<U_alpha,x^2,-1>}
\nonumber
\langle U_{\alpha}, x^{2,0} \rangle = \\ \nonumber
\left\langle \; \;  \left( \; \left( 
\begin{array}{c}
\alpha \cos \theta \cos t - \sin \theta \sin t \\
0   \end{array} \right) \right. \right.  , \; 
\left( \left. \begin{array}{c}
- \cos  \theta \sin t + \alpha \sin \theta \cos t \\
0  
\end{array} \right) \; \right) , \\ \nonumber
 \left( \; \left(\begin{array}{c} 
- \cos \theta \sin t  + \alpha \sin \theta \cos t \\
\sin \theta \sin t + \alpha \cos \theta \cos t
\end{array} \right) \right .   , \; 
\left. \left. \left( \begin{array}{c} 
\alpha \cos \theta \cos t - \sin \theta \sin t    \\
- \alpha \sin \theta \cos t - \cos \theta \sin t 
\end{array}\right) \; \right) \; \; \right\rangle        = \\
\frac{1}{2} ( \ \cos^2 t \sin \theta \cos \theta + \sin^2 t \cos \theta \sin \theta  + \sin^2 t \cos \theta \sin \theta  +
\cos^2 t \sin \theta \cos \theta \ ) =\nonumber \\ 
\frac{1}{2}  \sin 2 \theta, \qquad\mbox{}
\end{eqnarray}\normalsize  
Similar computations show that 
\addtocounter{thrm}{1} 
\begin{eqnarray}\label{projections of x and y onto U_gamma_i and D_gamma_i}
 \langle U_{\gamma_i}, x^{2,0} \rangle = \langle U_{\alpha},  \eta_{1}^{2,0} \rangle =
\langle U_{\alpha}, \vartheta_{1}^{2,0} \rangle =   0.
\end{eqnarray}  
Combining  (\ref{< (- vartheta, vartheta), U_gamma_1 >}),
(\ref{<U_alpha,x^2,-1>}), and   
(\ref{projections of x and y onto U_gamma_i and D_gamma_i}) 
we see that the projection of $P$ onto $TO_{A^u}$ is nondegenerate     
unless $\sin2 \theta = 0$ or $\cos 2 \theta = 0$, proving the proposition.
\end{proof}

\noindent
\rk{Remark}
One might hope that there really is positive curvature on these planes when $\theta = 0,
\frac{\pi}{4}, \frac{\pi}{2}$ and $\frac{3\pi}{4}$ because perhaps the projection 
of these planes onto $TO_{A^d}$ is nondegenerate.  While this is true 
for most values of $\theta$, it is false precisely when $\theta = 0,
\frac{\pi}{4}, \frac{\pi}{2}$ and $\frac{3\pi}{4}$.  \vspace{.1in}

Finally we study the zero curvatures of $(E_{2,0}, g_{\nu_1, \nu_2, l_{1}^u, l_{2}^d})$
at points when $t = \frac{\pi}{4}$.   At such points the action induced by $A_{SO(2)}$
on $S^4$ via $p_{2,0}\co E_{2,0} \longrightarrow S^4$ fixes $S^{2}_{im}$, so the  
orbits of $A_{SO(2)}$ are vertical with respect to $p_{2,0}$
 when $t= \frac{\pi}{4}$.  Combining this with the fact that $A^l$ acts 
transitively on the fibers of $p_{2,0}$,  we see that 
to study the curvature of these points we may assume that $\theta=0$,
 in other words we are back to representative points of the form 
(\ref{modulo isometries}) with $t=\frac{\pi}{4}$.  
At such points we get the following as a
corollary of (\ref {<(0, vartheta), U_gamma_1>}),  
(\ref{ < U_gamma_1, (vartheta_1, 0) >}), 
(\ref { < U_gamma_2, vartheeta_1 >}) and  
(\ref { < U_gamma_2, (vartheta_1, 0) >}).

\begin{coro}\label{t = pi/4}
When $t = \frac{\pi}{4}$ and $\theta = 0$  
the planes of the form (\ref{extra 0's when t = pi/4}) have positive curvature
except if they have the form
\addtocounter{thrm}{1} 
\begin{eqnarray}\label{special vartheta 0's} 
P = span\{ ( \ \vartheta, 0 \ ), \ ( 0, \vartheta ) \}
\end{eqnarray}
where $\vartheta$ is any convex combination of $\{ \vartheta_1, \vartheta_2 \}$. 
\end{coro}

To determine which of the planes in (\ref{x-y zeros at pi/4}) still have $0$--curvature
we prove the following.
\begin{prop}\label{stable x-y 0's at pi/4}
At points with $t= \frac{\pi}{4}$ and $\theta= 0$, the only planes of the form
(\ref{x-y zeros at pi/4}) that have $0$--curvature with respect to
$g_{\nu_1, \nu_2, l_{1}^u, l_{1}^d}$ are those of the form
\addtocounter{thrm}{1}
\begin{eqnarray}\label{stable x o's at pi/4}
P= span\{ x^{2,0} + \lambda ( - \vartheta, 0) , ( - \vartheta, \vartheta) \}   
\end{eqnarray}
and
\addtocounter{thrm}{1}
\begin{eqnarray}\label{stable y o's at pi/4}
P= span\{ y^{2,0} + \lambda ( - \vartheta, 0) , ( \vartheta, \vartheta) \}   
\end{eqnarray}
where $\vartheta$ stands for any convex combination of $\{ \vartheta_1, \vartheta_2 \}$
and $\lambda$ is any real number.  
\end{prop}

\begin{proof}
The point is that if $z$ is any other convex combination of $\{ x^{2,0}, y^{2,0} \}$,
then the projection of the plane in (\ref{x-y zeros at pi/4}) onto $TO_{A^u}$ is 
nondegenerate.  To see this we will need to compute the projection of 
$y^{2,0}$ onto $TO_{A^u}$.  Although we will only need to know its values
for $( t, \theta) = ( \frac{\pi}{4}, 0)$, we will do the computation for 
arbitrary $(t, \theta)$, since it may be of interest to some readers.  

Along the orbit of $SO(2)$ the formula for $y^{2,0}$ is \small
\begin{eqnarray*}
\left(\begin{array}{cc}
\cos \theta    &    \sin \theta   \\
- \sin \theta  &    \cos \theta   
\end{array}\right) 
\left( \; \left(\begin{array}{c} 
	\alpha \sin t    \\ 
\cos t    \end{array} \right) \  \right. , \;  
\left.   \left( \begin{array}{c}   
	- \cos t   \\
	- \alpha \sin t   \end{array} \right) \; \right) =  \\   
\left( \; \left(\begin{array}{c} 
\alpha \cos \theta \sin t  +  \sin \theta \cos t \\  
- \alpha \sin \theta \sin t +  \cos \theta \cos t
\end{array} \right) \right .   \ , \; 
\left. \left( \begin{array}{c}   
- \cos \theta \cos t - \alpha \sin \theta \sin t    \\
\sin \theta \cos t - \alpha \cos \theta \sin t 
\end{array}\right) \; \right).        
\end{eqnarray*}\normalsize
So 
\addtocounter{thrm}{1}\small
\begin{eqnarray}\label{<U_alpha, y>}  
\nonumber
\langle U_{\alpha}, y^{2,0} \rangle =  \\
\nonumber
\left\langle \; \;  \left( \; \left( 
\begin{array}{c}
\alpha \cos \theta \cos t - \sin \theta \sin t \\
0   \end{array} \right) \right. \right. \  , \; 
\left( \begin{array}{c} 
- \cos  \theta \sin t + \alpha \sin \theta \cos t \\
0    
\end{array} \right) , \\
\nonumber
\left( \; \left(\begin{array}{c} 
\alpha \cos \theta \sin t  +  \sin \theta \cos t \\
- \alpha \sin \theta \sin t +  \cos \theta \cos t
\end{array} \right) \right .   \ , \;   
\left. \left. \left( \begin{array}{c} 
- \cos \theta \cos t - \alpha \sin \theta \sin t    \\
\sin \theta \cos t - \alpha \cos \theta \sin t 
\end{array}\right) \; \right) \; \; \right\rangle = \\
\frac{1}{2} ( \  \cos^2 \theta \sin t \cos t - \sin^2 \theta \cos t \sin t + 
\cos^2 \theta \sin t \cos t  - \sin^2 \theta \cos t \sin t \ )   = 
\nonumber\\   
\frac{1}{2} \sin t \cos t ( 2 \cos^2 \theta - 2 \sin^2 \theta ) =  
\nonumber  \\
\frac{1}{2}  \sin 2 t \cos 2 \theta.  \qquad\mbox{}
\end{eqnarray}\normalsize
 Similar computations show that 
\addtocounter{thrm}{1}
\begin{eqnarray}\label{projection of y onto U_gamma_i}
 \langle U_{\gamma_i}, y^{2,0} \rangle = \langle U_{\alpha}, ( \vartheta, 0 )\rangle =
\langle U_{\alpha}, (0, \vartheta) \rangle =  0.
\end{eqnarray}   
Combining 
(\ref{<U_alpha, y>}) and (\ref{projection of y onto U_gamma_i}) with
(\ref{<(0, vartheta), U_gamma_1>}), (\ref{ < U_gamma_1, (vartheta_1, 0) >}), 
(\ref{ < U_gamma_2, vartheeta_1 >}), (\ref{ < U_gamma_2, (vartheta_1, 0) >}),  
(\ref{<U_alpha,x^2,-1>}) and 
(\ref{projections of x and y onto U_gamma_i and D_gamma_i}) we see that 
if $P$ is any plane of the form of (\ref{x-y zeros at pi/4}) and not of  
the form of either (\ref{stable x o's at pi/4}) or
(\ref{stable y o's at pi/4}), then the projection of $P$ onto 
$TO_{A^u}$ is nondegenerate and hence $P$ is positively curved.  
\end{proof} 
  
\noindent
\rk{Remark}
As was the case with our previous results one could study projections
onto $TO_{A^d}$ to prove (\ref{t = pi/4}) and (\ref{stable x-y 0's at pi/4}).  Unfortunately
these projections are degenerate on the same planes as the projections 
onto $TO_{A^u}$ and hence do not give us additional positive curvature.   
\vspace{.1in}

The curvature computations for the proof of Theorem \ref{unit tangent bundle}
are completed with the following.

\begin{prop}\label{really are 0's}
\mbox{}\newline
\vspace{-.18in} \mbox{}
\begin{itemize}
\item[\bf(\romannumeral1)]
The planes in (\ref{generic eta-vartheta 0's}), (\ref{generic x, vartheta 0's}), 
(\ref{special vartheta 0's}), (\ref{stable x o's at pi/4}) and 
(\ref{stable y o's at pi/4}) all have $0$--curvature in      
$(E_{2,0}, g_{\nu_1, \nu_2, l_{1}^u, l_{1}^d} )$.   
\item[\bf(\romannumeral2)]
The plane spanned by $x^{2,0}$ and $ \vartheta_{1}^{2,0}$ in $E_{2,0}$
is tangent to a totally geodesic flat $2$--torus.  
\item[\bf(\romannumeral3)] 
	Let $Z$ be the pointwise zero locus in $(E_{2,0}, g_{\nu_1, \nu_2, 
l_{1}^u, l_{1}^d } )$, then $Z$ is the union of two copies of 
$S^3 \times S^3$ that intersect along a copy of $S^2 \times S^3$.  
\end{itemize}
\end{prop}

\begin{proof}[Sketch of proof]
To prove (\romannumeral1) we must check that the $A$--tensors of 
$q_{A^u \times A^d}^{diag}$ and $q_{2,0}$ vanish on these    
planes.  These computations are rather long, but are quite 
straight forward so we leave them to the reader.  

Recall \cite{Wil1} that $S^{1}_{\Bbb R}$ denotes the circle in $S^4$ 
that is fixed by the action induced from $A^r$ via $p_{2,0}$.

It is easy to see that $x^{2,0}$ and $ \vartheta_{1}^{2,0}$
span a $2$--torus, $T$, in $Sp(2)$ whose tangent plane is horizontal with respect to 
$q_{2,0}$ at every point.  Since $T$ is horizontal, $q_{2,0}|_{T}$ is    
a covering map.   The order of the covering is at least two because 
$-1$ leaves $T$ invariant under $A_{2,0}$.  We show next that the order is 
exactly two.  

 Let $c_x$ and $c_{\vartheta_{1}^{2,0}}$ be the geodesic circles 
that are generated by $x^{2,0}$ and $ \vartheta_{1}^{2,0}$.   Then with respect to 
the join decomposition $S^4 = S^{2}_{im} * S^{1}_{\Bbb R}$, $p_{2,0}(c_x)$  
is a radial circle and for $t> 0$, $p_{2,0}(c_{\vartheta_{1}^{2,0}})$ is 
an intrinsic geodesic in a copy of $S^2$.  Using these observations we see 
that a point in $p_{2,0}( q_{2,0}( T))$ has only two preimages under 
$p_{2,0} \circ q_{2,0}$ from which it follows that $q_{2,0}|_T$ is 
a two fold covering.    

The fields $x^{2,0}$ and $\vartheta_{1}^{2,0}$ are both parallel along 
$T$.  Straightforward computation shows that the results of 
transporting them from a point $z \in T$ to $-z$ agrees with 
their images under the differential of the map 
induced on $T$ by $-1$ via $A_{2,0}$.  Therefore $q_{2,0}|_T$ 
is a 2--fold orientation preserving cover, and $q_{2,0}(T)$ is a 
torus rather than a Klein bottle.

It follows from (\romannumeral1) that $q_{2,0}(T)$ is  flat and computations similar to those 
in the 
proof of (\romannumeral1) show that $q_{2,0}(T)$ is totally geodesic.

Combining part (\romannumeral1) with (\ref{0's wrt g_nu_1, nu_2}), 
(\ref{proj onto  eta}) and (\ref{x-y, planes}), we see that $Z$ consists 
of the points with
$\theta =0, \frac{\pi}{4}, \frac{\pi}{2},$ or $\frac{3\pi}{4}$ or $t = \frac{\pi}{4}$.  
Using the join structure $S^4 = S^{1}_{\Bbb R} * S^{2}_{im}$ 
we see that the union of the points with $\theta = 0$,
$\theta = \frac{\pi}{2}$ and $t= \frac{\pi}{4}$ is a $3$--sphere in $S^4$.   
(Keep in mind that $A_{SO(2)}$ is ${\Bbb Z}_2$--ineffective on 
$S^4$.)  For similar reasons  the union of the points with $\theta = \frac{\pi}{4}$, 
$\theta = \frac{3\pi}{4}$ and $t= \frac{\pi}{4}$ is another $3$--sphere.  These 
two three spheres intersect along the common $2$--sphere $S^{2}_{im}$
and their inverse images via $p_{2,0}$ are of course diffeomorphic to
$S^3 \times S^3$.  Finally notice that every point in this 
inverse image has zero curvatures because $p_{2,0}$ is the  
quotient map of $A^r$ and $A^r$ acts by isometries.   
\end{proof}

The example, $M^6$, of Corollary \ref{6-dim quot} is 
$$
M^6 = E_{2,0}/A^{l}_{\beta}
$$
where $\beta$ is any purely imaginary, unit quaternion, and $A^{l}_{\beta}$ 
is the circle subaction of $A^l$ that is generated by $\beta$.  
When $\alpha = \beta$ the torus $q_{2,0}(T) \subset E_{2,0}$ 
is normal to the orbit of $A^{l}_{\beta}$.  Some further curvature 
computations show that its image in $M^6$ is totally geodesic and flat.

\section{Topological Computations}

\begin{prop}\label{E_2,0 is the unit tangent bunlde}
$E_{2,0}$ is diffeomorphic to the total space of the unit
tangent bundle of $S^4$.  In fact, the Riemannian submersion
$p_{1,1} \co  E_{2,0} \longrightarrow S^4$ that is given by 
\begin{eqnarray*}
p_{1,1}\co  orbit \left( \ \left( \begin{array}{cc}
		a  &   b   \\
		c  &   d   \end{array}\right) \ \right)
\mapsto \tilde{h} ( c, d)  
\end{eqnarray*}
is bundle isomorphic to the unit tangent bundle of $S^4$, where 
$\tilde{h}\co  S^7 \longrightarrow S^4$ is the Hopf fibration given by 
left multiplication.   
\end{prop}

\begin{proof}
Translating Theorem 9.5 on page 99 of \cite{Huse} into 
our classification scheme (\ref{gluing map}) shows 
that the unit tangent bundle is the bundle of type 
$(1,1)$.  We will show via direct computations (similar to
those in \cite{GromMey}) that $( E_{2,0}, p_{1,1})$ is also the 
bundle of type $(1,1)$.  

As in \cite{Wil1} we define $\phi \co  {\Bbb R}^4 \longrightarrow {\Bbb R}$ 
by 
$$
\phi(u) = \frac{1}{ \sqrt{ 1 + |u|^2} },    
$$
and we define explicit bundle charts 
$h_1, h_2 \co  {\Bbb R}^4 \times S^3 \longrightarrow E_{2,0}$
by 
\begin{eqnarray*}
h_1 (u,q) = orbit \left( \begin{array}{cc}  
		- q          &      qu    \\ 
		\bar{u}    &       1    \end{array} \right)  \phi(u)  
\end{eqnarray*}
and
\begin{eqnarray*}
h_2 (v, r) = orbit \left( \begin{array}{cc} 
		  - r \bar{v}   &     r     \\
		 1              &    v    
				\end{array} \right) \phi(v).
\end{eqnarray*}
$h_1$ and $h_2$ are embeddings onto the open sets 
\begin{eqnarray*}
U_1 = orbit \{ \ \left( \begin{array}{cc} 
			a   &   b    \\
			c   &   d   \end{array} \right)  \; | \;
		d \not= 0   \; \}    
\end{eqnarray*}
and 
\begin{eqnarray*}
U_2 = orbit \{  \ \left( \begin{array}{cc} 
			a   &   b    \\ 
			c   &   d   \end{array} \right)  \; | \;
		c \not= 0   \; \}.   
\end{eqnarray*}
In fact their inverses are given by 
\begin{eqnarray*}
h_{1}^{-1} ( \ orbit \left( \begin{array}{cc} 
			a  &   b   \\
			c  &    d   \end{array} \right)  \ ) = 
( \ \frac{\bar{c} d} { | d |^2 } , \ \frac{ - \bar{d} a }{ | d| | a| } \ )    
\end{eqnarray*}
and 
\begin{eqnarray*}
h_{2}^{-1} ( \ orbit \left( \begin{array}{cc} 
			a  &   b   \\ 
			c  &    d   \end{array} \right)  \ ) = 
( \ \frac{\bar{c} d} { | c |^2 } , \ \frac{\bar{c} b }{ | c| | b| } \ ).
\end{eqnarray*}
Thus
\begin{eqnarray*}  
h_{2}^{-1} \circ h_1 ( u, q) = 
h_{2}^{-1} (   orbit \left( \begin{array}{cc}   
		  - q           &     q u    \\   
		 \bar{u}        &      1    
				\end{array} \right) \phi(u) )=  \\
(  \  \frac{u}{|u|^2}   , \ \frac{u q u}{| u|^2 }  \ ).
\end{eqnarray*} 
So $(E_{2,0}, p_{1,1})$ is the bundle of type $(1,1)$ and hence 
is the unit tangent bundle of $S^4$. 
\end{proof}

Next we compute the homology and first few homotopy groups of $S^3$--bundles
over $S^4$ to distinguish $E_{2,0}$ from the known examples 
of $7$--manifolds with positive sectional curvature.
\begin{prop}\label{homology of S^3-bundles over S^4}
Let $E_{m,-n}$ denote the $S^3$--bundle over $S^4$ of type $(m,-n)$.  
\begin{itemize}
\item[\bf(\romannumeral1)]
		The integral homology groups of $E_{m,-n}$ are 
\begin{eqnarray*}
 H_0(E_{m,-n}, {\Bbb Z}) \cong H_7(E_{m,-n}, {\Bbb Z}) \cong {\Bbb Z} \\ 
H_3( E_{m,-n} , {\Bbb Z}) \cong \frac{\Bbb Z}{ (m-n) Z}, \; \; \mbox{and} \\
H_q(E_{m,-n}, {\Bbb Z}) \cong \{ 0\} \; \; \mbox{for all $q \not= 0, 3, 7$,}
\end{eqnarray*}
if $m \not= n$.  
\item[\bf(\romannumeral2)]
$\pi_1(  E_{m,-n} ) \cong \pi_2 ( E_{m,-n} ) \cong \{ 0 \} \; \; \mbox{and} \; \; 
\pi_3( E_{m,-n} ) \cong  \frac{\Bbb Z}{ (m-n) Z}$
\item[\bf(\romannumeral3)]
$E_{2,0}$ does not have the homotopy type of any known example of 
a manifold with positive curvature.  
\end{itemize}
\end{prop}

\begin{proof}
To compute the homology of $E_{m,-n}$ we decompose it as 
\addtocounter{thrm}{1}  
\begin{eqnarray}\label{decomposition of E_m,-n}  
E_{m,-n} = D^4 \times S^3 \cup_{g_{m,-n}} D^4 \times S^3
\end{eqnarray}
where  the gluing map 
$g_{m,-n} \co  S^3 \times S^3 \longrightarrow S^3 \times S^3$
is given by $g_{m,-n}(u, v) = (u, u^m v u^{-n} )$.  
From (\ref{decomposition of E_m,-n}) and the 
Seifert--Van Kampen theorem it follows  
that $E_{m,-n}$ is simply connected.   

The 
Mayer--Vietoris sequence for the decomposition 
(\ref{decomposition of E_m,-n})
is
\addtocounter{thrm}{1}
\begin{eqnarray}\label{Meyer-Vietoris}
\nonumber
\cdots \rightarrow H_{q}( S^3 \times S^3)
\stackrel{\Phi_q }{\rightarrow} H_q( D^4 \times S^3) \oplus 
 H_q( D^4 \times S^3)\\ \stackrel{ \Psi_{q} }{\rightarrow }
H_q(E_{m,-n}) \stackrel{ \Gamma_q}{\rightarrow} H_{q -1}(
S^3 \times S^3 ) \stackrel{\Phi_{q-1}}{\rightarrow} \cdots
\end{eqnarray}
Since $H_2 (D^4 \times S^3 ) \cong H_1( S^3 \times S^3) \cong 
\{ 0 \}$, $H_2( E_{m,-n} ) \cong \{ 0\}$.  

Also since $H_2( S^3 \times S^3) \cong \{ 0 \}$, 
$H_3( E_{m,-n} ) \cong 
\frac{H_{3}(D^4 \times S^3) \oplus H_{3}(D^4 \times S^3) }
{ image(\Phi_3) }$.  The next step is to compute $image(\Phi_3)$.  

Suppose $\{ \alpha \times 1, 1 \times \alpha \}$ is a set of 
generators for $H_3 (S^3 \times S^3)$, and $\{ ( \beta, 0 ), (0, \beta)\}$
generates $H_3( D^4 \times S^3) \oplus  H_3( D^4 \times S^3)$.  
Then using the fact that the map $P_k \co  S^3 \longrightarrow S^3$ given by
$P_k( q) = q^k$ has degree $k$ we can see that 
\begin{eqnarray*}
\begin{array}{ll}
\Phi_3 ( \alpha \times 1 ) = (m -n ) (0, \beta)   & \mbox{and} \\
\Phi_3 ( 1 \times \alpha ) = ( \beta, 0 ) +  (0, \beta).    & \mbox{} 
\end{array}
\end{eqnarray*}
( To evaluate the first map represent $\alpha \times 1$ by  
$S^{3} \times \{ 1 \}$.  )

Put another way, the matrix of $\Phi_3$ with respect to our sets 
of generators is 
\addtocounter{thrm}{1}
\begin{eqnarray}\label{matrix of Phi_3}
\left(\begin{array}{cc}
0    &    1     \\
m -n     &    1
\end{array}\right). 
\end{eqnarray}    
From this point a routine algebraic computation shows that 
$$
H_3( E_{m,-n} ) \cong 
\frac{H_{3}(D^4 \times S^3) \oplus H_{3}(D^4 \times S^3) }
{ image(\Phi_3) } \cong \frac{ {\Bbb Z} }{ (m -n) {\Bbb Z} }.  
$$
It follows from (\ref{matrix of Phi_3}) that $\Phi_3$
is injective if $m \not= n$, and hence we get using (\ref{Meyer-Vietoris})
that 
$$
H_4( E_{m,-n} ) \cong H_5( E_{m,-n} ) \cong H_6( E_{m,-n} ) \cong
0,
$$
completing the proof of (\romannumeral1).  

Part (\romannumeral2) is a corollary of (\romannumeral1). 

The known examples of simply connected $7$--manifolds with positive curvature are 
given in \cite{Berger1}, \cite{Esch1} and \cite{AllWal}.   With 
the exception of $S^7$ and the example, $V_1$, of Berger
none of them are $2$--connected.   Since $H_3( E_{2,0} , {\Bbb Z}) \cong { \Bbb Z}_2$,
$E_{2,0}$ is not a homotopy sphere.  It is 
not homotopy equivalent $V_1$
since, according to Proposition 40.1 in \cite{Berger1},  
  $V_1$ is not a ${\Bbb Z}_5$--cohomology sphere.   
\end{proof}    
    
Our last topological computation is the following (cf also
\cite{GrovZil} Corollary 3.9).
\begin{prop}\label{cohomology of M^6}
\mbox{}\newline
\vspace{-.18in} \mbox{}
\begin{itemize}
\item[\bf(\romannumeral1)]
$M^6$ has the same integral cohomology modules as $CP^3$ 
 but not the same integral cohomology algebra.
\item[\bf(\romannumeral2)]
$M^6$ does not have the homotopy type of any known example 
of a manifold of positive curvature.
\end{itemize}
\end{prop}

\begin{proof}
From the long exact homotopy sequence of the fibration $%
S^{1}\hookrightarrow E_{2,0}\longrightarrow M^{6}$ we see that $\pi
_{1}(M^{6})\cong 0$ and $\pi _{2}(M^{6})\cong \Bbb{Z}$. Thus $H^{1}(M^{6},%
\Bbb{Z})\cong 0$ and $H^{2}(M^{6},\Bbb{Z})\cong \Bbb{Z}$. To compute $%
H^{3}(M^{6},\Bbb{Z})$ and the cup products we appeal to the Gysin
sequence 
\[
\cdots \rightarrow H^{1}(M^{6},\Bbb{Z})\rightarrow H^{3}(M^{6},\Bbb{Z}%
)\rightarrow H^{3}(E_{2,0},\Bbb{Z})\rightarrow H^{2}(M^{6},\Bbb{Z}%
)\rightarrow \cdots 
\]
of the fiber bundle $S^{1}\hookrightarrow E_{2,0}\longrightarrow M^{6}$ with
integer coefficients. We already know that $H^{1}(M^{6},\Bbb{Z})=0$.
Moreover, it follows from the universal coefficient theorem that $%
H^{3}(E_{2,0},\Bbb{Z})=0$ since $H_{2}(E_{2,0},\Bbb{Z})$ has no torsion and $%
H_{3}(E_{2,0},\Bbb{Z})$ is a torsion group. Thus the Gysin sequence shows
that $H^{3}(M^{6},\Bbb{Z})=0$ and hence that $M^{6}$ has the same integral
cohomology groups as $\Bbb{C}P^{3}.$ To see that $M^{6}$ does not have the
cohomology ring of $\Bbb{C}P^{3}$ we look at the Gysin sequence again: 
\begin{eqnarray*}
0=H^{3}(E_{2,0},\Bbb{Z})\stackrel{\pi _{\ast }}{\rightarrow }H^{2}(M^{6},%
\Bbb{Z})\stackrel{\cup e}{\rightarrow }H^{4}(M^{6},\Bbb{Z})\stackrel{\pi
^{\ast }}{\rightarrow }\hskip 1.085in\mbox{}\\
H^{4}(E_{2,0},\Bbb{Z})\stackrel{\pi _{\ast }}{%
\rightarrow }H^{3}(M^{6},\Bbb{Z})=0,\\
0=H^{5}(E_{2,0},\Bbb{Z})\stackrel{\pi _{\ast }}{\rightarrow }H^{4}(M^{6},%
\Bbb{Z})\stackrel{\cup e}{\rightarrow }H^{6}(M^{6},\Bbb{Z})\stackrel{\pi  
^{\ast }}{\rightarrow }H^{6}(E_{2,0},\Bbb{Z})=0.
\end{eqnarray*}
We therefore have 
\begin{eqnarray*}
&0 \rightarrow H^{2}(M^{6},\Bbb{Z})\stackrel{\cup e}{\rightarrow }%
H^{4}(M^{6},\Bbb{Z})\stackrel{\pi ^{\ast }}{\rightarrow }\Bbb{Z}%
_{2}\rightarrow 0, \\
&0 \rightarrow H^{4}(M^{6},\Bbb{Z})\stackrel{\cup e}{\rightarrow }%
H^{6}(M^{6},\Bbb{Z})\rightarrow 0.
\end{eqnarray*}
This means that $H^{2}(M^{6},\Bbb{Z})\stackrel{\cup e}{\rightarrow }%
H^{4}(M^{6},\Bbb{Z})$ maps the generator $x$ for $H^{2}(M^{6},\Bbb{Z})$ to
twice a generator $y$ for $H^{4}(M^{6},\Bbb{Z}).$ If we let $e=kx$ we
therefore have $x\cup kx=2y.$ Therefore, it suffices to show that $k=\pm 1$
in order to show that the cohomology ring of $M^{6}$ is not that of $\Bbb{C}%
P^{3}.$ We have that $H^{4}(M^{6},\Bbb{Z})\stackrel{\cup e}{\rightarrow }%
H^{6}(M^{6},\Bbb{Z})$ is an isomorphism. Thus $( x\cup kx ) \cup kx=k^{2}x\cup
x\cup x$ must also be twice a generator for $H^{6}(M^{6},\Bbb{Z}).$
Therefore, if $k=\pm 2$ we would have that $\pm 2x\cup x\cup x$ generates $%
H^{6}(M^{6},\Bbb{Z}).$ This, however, is impossible.

Besides $S^{6}$ and $\Bbb{C}P^{3}$, there are two examples given in \cite{Wal}
and \cite{Esch2} of simply connected $6$--manifolds with positive curvature.
Neither of these has the cohomology modules of $\Bbb{C}P^{3}$ so $M^{6}$ does not
have the homotopy type of any known example. 
\end{proof}

\end{document}